\documentclass{amsart}

\usepackage{amsfonts, amsmath, amssymb, amsthm, amsxtra, latexsym}

\theoremstyle{plain}

\newtheorem{theorem}{Theorem}[section]
\newtheorem{proposition}[theorem]{Proposition}
\newtheorem{lemma}[theorem]{Lemma}
\newtheorem{corollary}[theorem]{Corollary}
\newtheorem{conjecture}[theorem]{Conjecture}

\theoremstyle{definition}

\newtheorem{remark}[theorem]{Remark}

\numberwithin{equation}{section}

\renewcommand{\Re}{{\textup{Re\,}}}
\renewcommand{\Im}{{\textup{Im\,}}}
\newcommand {\R}{{\mathbb R}}
\newcommand {\C}{{\mathbb C}}
\newcommand{\Z}{{\mathbb Z}}

\newcommand{\ga}{{\mathfrak a}}
\newcommand{\gac}{\ga_\C}
\newcommand{\gp}{{\mathfrak p}}
\newcommand{\gq}{{\mathfrak q}}

\newcommand{\e}{\epsilon}
\newcommand{\g}{\gamma}
\newcommand{\G}{\Gamma}
\newcommand{\D}{\Delta}
\newcommand{\x}{\xi}
\newcommand{\n}{\nu}
\newcommand{\p}{\pi}

\newcommand{\tu}{\textup}
\newcommand{\on}{{{\textup O}(N)}}

\newcommand{\poly}{{\mathcal P}}
\newcommand{\schw}{{\mathcal S}}
\newcommand{\test}{{\mathcal D}}

\newcommand{\reg}{{\textup{reg}}}
\newcommand{\sing}{{\textup{sing}}}
\newcommand{\ksing}{K^{\sing}}
\newcommand{\kreg}{K^{\reg}}
\newcommand{\gareg}{{\ga^\reg}}
\newcommand{\gacreg}{{\gac^\reg}}
\newcommand{\refl}{r_\alpha}

\newcommand{\rplus}{{R_+}}
\newcommand{\sumposroots}{\sum_{\alpha\in\rplus}}
\newcommand{\tk}{T_\x(k)}

\newcommand{\Exp}{{\textup{Exp}_G}}
\newcommand{\w}{w_k}

\newcommand{\inta}{\int_{\ga}}
\newcommand{\intp}{\int_{\gp}}
\newcommand{\intq}{\int_{\gq}}

\newcommand{\dk}{D_k}
\newcommand{\ek}{E_k}

\newcommand{\lone}{L_1(\ga,\vert \w (x)\vert dx)}

\newcommand{\mxi}{M_{\x^*}}
\newcommand{\deltak}{\D_k}
\newcommand{\expdeltak}{e^{\frac{\deltak}{2}}}
\newcommand{\expdeltazero}{e^{\frac{\D_0}{2}}}
\newcommand{\expmindeltak}{e^{-\frac{\deltak}{2}}}
\newcommand{\expmindeltakp}{e^{-\frac{\D_{k^\prime}}{2}}}
\newcommand{\Aut}{{\textup{Aut}}}
\newcommand{\End}{{\textup{End}}}
\newcommand{\Hom}{{\textup{Hom}}}
\newcommand{\Span}{{\textup{Span}}}
\newcommand{\Res}{{\textup{Res}\,}}
\newcommand{\Resxo}{{\textup{Res}_{x_0}\,}}
\newcommand{\Ext}{{\textup{Ext\,}}}

\newcommand{\ad}{{\textup{ad}}}
\newcommand{\ext}{{\textup{ext}}}
\newcommand{\supp}{{\textup{supp}}}
\newcommand{\Rad}{{\textup{Rad}}}

\newcommand{\mpsi}{M_\psi}

\newcommand{\ball}{{B_R}}
\newcommand{\openball}{B_R^\prime}
\newcommand{\testb}{\test(\ball)}
\newcommand{\testbdense}{\test_0(\ball)}
\newcommand{\testbinv}{\testb^{\on}}
\newcommand{\schwinv}{\schw(\ga)^\on}
\newcommand{\schwdense}{\schw_0(\ga)}
\newcommand{\schwinvdense}{\schw_0(\ga)^\on}
\newcommand{\pwtype}{{\mathcal H}_\ball}
\newcommand{\harm}{{\mathcal H}}

\newcommand{\dist}{\test^\prime}
\newcommand{\tempdist}{\schw^\prime}
\newcommand{\onehalf}{\frac{1}{2}}

\newcommand{\group}{{\mathcal G}}
\newcommand{\maxcompact}{{\mathcal K}}
\newcommand{\restoa}{\Res_{\mathfrak a}^{\mathfrak p}}
\newcommand{\exttop}{\Ext_{\mathfrak a}^{\mathfrak p}}
\newcommand{\fourier}{{\mathcal F}}

\newcommand{\khat}{\widehat{K}}
\newcommand{\ondual}{\widehat O (N)}

\begin{document}


\title{Paley--Wiener theorems for the Dunkl transform}


\translator{}

\dedicatory{}

\author[Marcel~de~Jeu]{Marcel~de~Jeu}

\begin{abstract}
We conjecture a geometrical form of the Paley--Wiener theorem for the Dunkl
transform and prove three instances thereof, by using a reduction to the one-dimensional even case, shift
operators, and a limit transition from Opdam's results for the graded Hecke
algebra, respectively. These Paley--Wiener theorems are used to extend Dunkl's
intertwining operator to arbitrary smooth functions.

Furthermore, the connection between Dunkl operators and the Cartan motion group
is established. It is shown how the algebra of radial parts of invariant
differential operators can be described explicitly in terms of Dunkl operators.
This description implies that the generalized Bessel functions coincide with
the spherical functions. In this context of the Cartan motion group, the
restriction of Dunkl's intertwining operator to the invariants can be
interpreted in terms of the Abel transform. We also show that, for certain
values of the multiplicities of the restricted roots, the Abel transform is
essentially inverted by a differential operator.
\end{abstract}

\date{}



\thanks{}

\address{M.F.E.~de~Jeu\\
         Mathematical Institute\\
         Leiden University\\
         P.O. Box 9512\\
     2300 RA Leiden\\
     The Netherlands}



\email{mdejeu@math.leidenuniv.nl}

\urladdr{}

\maketitle

\section{Introduction and overview}\label{sec:intro}

In recent times the study of special functions associated with root systems has
developed to a considerable degree. Starting with a number of conjectures by
Macdonald, and the work of Heckman and Opdam on multivariable hypergeometric
functions in the late 1980's, the development of the theory was greatly
enhanced by the introduction of rational Dunkl operators by Dunkl
\cite{diffdiff}. Through various intermediate steps of generalization, these
operators can even be said to have ultimately provided crucial building blocks
for Cherednik's work on double affine Hecke algebras and the $q$-Macdonald
conjectures.

Originally, before the introduction of Dunkl operators, the idea when studying
special functions related to root systems was to consider root multiplicities
in the theory of spherical functions on Lie groups as parameters, and then to
develop a theory for Weyl group invariant objects in this more general
situation, without the aid of the presence of the group. It was this point of
view which underlay the Macdonald conjectures and which led Heckman and Opdam
to the development of their theory of hypergeometric functions in higher
dimension. One of the main technical problems in this context is the
description of the generalized radial parts of invariant differential
operators. Apart from an explicit formula for the generalized radial part of
the Laplacian---an expression which was in fact the starting point for Heckman
and Opdam---the other operators remain somewhat intangible.

This problem disappeared when Dunkl found the operators which have come to bear
his name. They are parameterized deformations of the ordinary derivatives,
involving a finite reflection group, for which it is still relatively easy to
study the spectral problem and develop the theory of the corresponding Fourier
transformation---the Dunkl transform. The invariant part of the theory then
answers questions in a generalized theory of spherical functions---for the
Cartan motion group, to be precise---as described above. The operators under
consideration being explicitly given, it is actually \emph{easier} to study the
general context of Dunkl operators, and then specialize to invariant objects
later on, than it is to restrict oneself to the invariants from the outset. The
same holds true for the modification of Dunkl operators into Cherednik
operators, which gives rise to a representation of the graded Hecke algebra
\cite{gradhecke} and the non-invariant ``envelope'' of the work of Heckman and
Opdam on hypergeometric functions.

Quite remarkably, the theory of spherical functions in analysis on Lie groups
can in a number of situations thus be regarded as only the invariant part of a
general non-invariant theory for Dunkl-type operators. At the time of writing,
it is unknown to the author whether there is an underlying reason for this
phenomenon.

In view of all this, these Dunkl operators and their modifications have
attracted considerable attention in various areas of mathematics and
mathematical physics during the last decade. To get an impression of their
influence on the development of special functions associated with root
systems---also in the general non-invariant context---we refer to, e.g.,
\cite{heckmanbourbaki}, \cite{opdamjapan}, \cite{dunklxu},
\cite{roeslersummerschool}. For their use in the study of integrable quantum
many body systems of Calogero--Moser--Sutherland type we refer to
\cite{diejenvinet} and the bibliography therein.

In this paper, we are mainly concerned with the further development of the
general theory of the Dunkl transform, notably with the Paley--Wiener theorem.
In addition, we describe the relation between Dunkl operators and the Cartan
motion group. We will now turn to an overview of the contents of this paper
with regard to these two subjects.

The first substantial results for the Dunkl transform, i.e., the Plancherel
theorem and inversion theorem, were obtained by Dunkl \cite{hankeltransform}
and the author \cite{dunkltransform}. Two Paley--Wiener theorems were
established in \cite[Chapter 3]{thesis} (unpublished), of which this paper is
an extension. Whereas the proof of the inversion theorem in
\cite{dunkltransform} is partly a formal argument, based on various symmetry
properties of the eigenfunctions of Dunkl's operators, it will become apparent
below that quite some more work is required in Paley--Wiener theory, due to the
lack of adequate asymptotic results in the spectral domain. We will establish three Paley--Wiener
theorems, each of these being a special case of a conjectured geometrically
more precise general Paley--Wiener theorem, which can be found below as
Conjecture~\ref{conj:pw}. Each version requires a different technique. The
first and most general version relies on a reduction to the one-dimensional
even case---where asymptotic results are available---and some non-trivial
results from general representation theory and from the representation theory
of the orthogonal group. The second version, which is proved for a discrete set
of parameters only, is established through shift operators. Finally, the third
version, which we prove for Weyl groups, follows from a limit transition from
results of Opdam \cite{gradhecke}. This limit transition has some interest in
itself, and for a discrete set of parameters its validity has also been
established by Ben Sa\"id and \O rsted, cf.~\cite{saidoersted}.

Turning to the Cartan motion group, for which the material is taken from
\cite{thesis}, we will show how the algebra of radial parts of invariant
differential operators can be described explicitly in terms of Dunkl operators
with suitable values of the parameters. From this description it is easy to see
that the invariant components of the eigenfunctions of the Dunkl operators are
precisely the spherical functions. We can then also describe Dunkl's
intertwiner operator, or rather its restriction to the invariants, in terms of
the Abel transform. For certain multiplicities, the results on shift operators
imply that this Abel transform is essentially inverted by a differential
operator.

The organization of this paper is as follows. In
Section~\ref{sec:notations_previous} we establish the necessary general
definitions and notations, and we recall some previous results.
Section~\ref{sec:formularium} contains a number of useful formulas. We do not
claim any originality in particular for these formulas, but we believe that the
proofs are considerably simpler than the ones in the existing literature. The
conjectured Paley--Wiener theorem can be found in Section~\ref{sec:pw}, as can
the proven three instances of it which were described above. These
Paley--Wiener theorems are then used in Section~\ref{sec:intertwining} to
extend Dunkl's intertwining operator to the smooth functions. To conclude with,
Section~\ref{sec:cmg} contains the details of the connection between Dunkl
operators and the Cartan motion group, including the results on the Abel
transform.

\section{Notations and previous results}\label{sec:notations_previous}

Let $\ga$ be a real vector space of finite dimension $N$ which is equipped with
an inner product $(\,.\,,\,.\,)$, inducing a Lebesgue measure $dx$ on $\ga$. We
let $\gac=\ga\otimes_\R \C$ denote the complexification of $\ga$ and we extend
the form $(\,.\,,\,.\,)$ to a bilinear form on $\gac$, again denoted by
$(\,.\,,\,.\,)$. Both $\ga$ and $\gac$ can be identified with their duals via
$(\,.\,,\,.\,)$; for $\x$ in $\ga$ or $\gac$ the corresponding linear
functional is then denoted by $\x^*$. Define the orthogonal group
$\on=O(\ga,(\,.\,,\,.\,))$. The norm $\vert\,.\,\vert$ which is induced on
$\ga$ by $(\,.\,,\,.\,)$ extends to an $\on$-invariant norm on $\gac$, also
denoted by $\vert\,.\,\vert$. There is a natural action of $\on$ on functions:
$$
(g\cdot f)(x)=f(g^{-1}x).
$$
Let $G\subset\on$ be a finite (real) reflection group with corresponding root
system $R$. We may and will assume that $(\alpha,\alpha)=2$ for all $\alpha\in
R$. If $\alpha\in R$, then $\refl$ is the orthogonal reflection in the
hyperplane perpendicular to $\alpha$. We choose and fix a positive system
$\rplus$ in $R$.

A function $k:R\mapsto\C$ is called a multiplicity function if $k$ is
$G$-invariant. We write $k\geq 0$ if all values of $k$ are non-negative, with
analogous notations $\Re k\geq 0$ and $k>0$. We let $K$ denote the vector space
of multiplicity functions on $R$.

Let $k\in K$ and $\x\in\ga$. Then the corresponding Dunkl operator $\tk$
is defined by
$$
\tk=\partial_\x + \sumposroots k_\alpha (\alpha,\x) M_{(\alpha^*)^{-1}} (1-\refl).
$$
Here $\partial_\x$ ($\x\in\ga$) is the (unnormalized) directional derivative
operator; we have used the notation $M_f$ for pointwise multiplication by a function $f$.

The definition of the operators is independent of the choice of the positive system. The $\tk$
leave the polynomials $\poly$ invariant, mapping the homogeneous polynomials $\poly_n$ of degree $n$ into $\poly_{n-1}$. Furthermore, the spaces $C^\infty(\ga)$ of smooth functions, $\test(\ga)$ of compactly supported smooth functions and $\schw(\ga)$ of rapidly decreasing smooth functions are also invariant \cite[Lemma
2.1]{dunkltransform}.

Quite remarkably, the $\tk$ form a commutative family: $\tk T_{\eta}(k)=
T_{\eta}(k)T_\x(k)$ for all $\x,\,\eta\in\ga$ and all $k\in K$, as was proved
by Dunkl \cite{diffdiff}, see \cite{singpol} for a different proof. As a
consequence of the commutativity the map $\x\mapsto\tk$ extends to $\poly$; for
$p\in\poly$ the corresponding operator is then denoted by $T_p(k)$. The
important $k$-Laplacian $T_{\vert\,.\,\vert^2}(k)$ is denoted by $\D_k$.

For $k\geq 0$, Dunkl \cite{operatorscommuting} has constructed a linear
isomorphism $V_k:\poly\mapsto\poly$, homogeneous of degree $0$, such that $V_k
1=1$ and $T_\xi(k)V_k=V_k\partial_\xi$ $(\xi\in\ga)$. An alternative approach
for more general $k$ can be found in \cite{singpol}. A concrete description of
this intertwiner operator $V_k$ is presently still unknown, with the exception
of the one-dimensional case \cite{operatorscommuting} and the case $A_2$
\cite{dunkla2}. Significant abstract results were obtained by R\"osler
\cite{roeslerpositivity}, who showed, amongst others, that $V_k$ is for $k\geq
0$ a positive operator which can be described in terms of measures.

Let $\lambda\in\gac$ and consider the simultaneous eigenfunction problem

\begin{equation}\label{eq:eigenfunctionproblem}
\tk f=(\lambda,\x)f\quad(\x\in\ga).
\end{equation}
This problem was studied first by Dunkl \cite{integralkernels} for $k\geq 0$;
later Opdam \cite{besseletc} treated the general case. One more definition is
needed to state the general result: a multiplicity function $k\in K$ is said to
be \emph{singular} if the simultaneous kernel of $\{\tk\}_{\x\in\ga}$ in
$\poly$ is non-trivial, i.e., if it properly contains the constants. We will
use the self-evident notations $\ksing$ and $\kreg$.

The set $\ksing$ has been determined in all cases \cite{singpol} and some
partial information about the simultaneous kernel for singular multiplicities can also
be found in [loc.cit.]. The general nature of the simultaneous kernel for singular multiplicities is, however, still unknown, with the exception of the case $A_n$ of the symmetric groups which has been solved by Dunkl \cite{dunklsymmetric1, dunklsymmetric2}. In this paper we will mainly be concerned with
multiplicities satisfying $\Re{k}\geq 0$. Such multiplicities are regular, as
is most easily seen by considering the operator $\sum_{i=1}^N e_i^*T_{e_i}(k)$
for an orthonormal basis $\{e_1,\dots, e_N\}$ of $\ga$ \cite{singpol}.

For regular multiplicities, the result from \cite{besseletc} for the eigenfunction
problem is as follows.

\begin{theorem} For all $k\in\kreg$
the eigenfunction problem \eqref{eq:eigenfunctionproblem} has a $1$-dimensional solution
space for all $\lambda\in\gac$. This space contains a (unique) function $\Exp(\lambda,k,\,.\,)$ such that
$\Exp(\lambda,k,0)=1$. Furthermore, $\Exp(\lambda,k,\,.\,)$ extends to a holomorphic
function on $\gac$, and
$$
\Exp(\,.\,,\,.\,,\,.\,):\gac\mapsto\C
$$
is a meromorphic function with poles precisely in $\ksing$.
\end{theorem}

In order to be able to define the corresponding Fourier transform---the Dunkl
transform---for sufficiently general functions, one needs non-trivial bounds
for the eigenfunctions. It is shown in \cite{dunkltransform} that for
$\Re{k}\geq 0$ one has
\begin{equation}\label{eq:growthestimates}\vert\Exp(\lambda,k,z)\vert
\leq\sqrt{\vert G\vert}\exp(\max_{g\in G}\Re({g\lambda,z)}) \quad
(\lambda,z\in\gac),
\end{equation}
in particular
\begin{equation}\label{eq:growthestimates2}
\vert\Exp(i\lambda,k,x)\vert\leq\sqrt{\vert G\vert}\quad(\lambda,x\in\ga).
\end{equation}
If $k\geq 0$ then the constant $\sqrt{\vert G\vert}$ in \eqref{eq:growthestimates} and \eqref{eq:growthestimates2} can in fact be improved to $1$, as a consequence of R\"osler's results on the intertwiner operator \cite{roeslerpositivity}.

Following Dunkl, we define for $\Re{k}\geq 0$ the $G$-invariant complex-valued weight
function $\w=\prod_{\alpha\in\rplus} \vert\alpha^*\vert^{2k_\alpha}$. The $\tk$
are anti-symmetric with respect to this weight function: if $f\in\schw(\ga)$
and $g$ is smooth such that both $g$ and $\tk g$ are of at most polynomial
growth, then
\begin{equation}\label{eq:antisymmetry}
\inta (\tk f)\,g \w\,dx=- \inta f\,(\tk g)\w\,dx.
\end{equation}

In view of \eqref{eq:growthestimates2}, the Dunkl transform for $\Re{k}\geq 0$ is meaningfully defined on $\lone$ as
$$
\dk f(\lambda)=\frac{1}{c_k}\inta f(x)\Exp(-i\lambda,k,x)\w(x)\,dx
\quad(\lambda\in\ga,\,\,f\in\lone).
$$
Its alleged inverse is
$$
\ek f(x)=\frac{1}{c_k}\inta f(\lambda)\Exp(i\lambda,k,x)\w(\lambda)\,d\lambda
\quad(x\in\ga,\,\,f\in L_1(\ga,\vert\w(\lambda)\vert\,d\lambda)).
$$
Here the normalizing constant $c_k$ is the Mehta integral, which is defined as
$$
c_k=\inta \psi\w\,dx.
$$
There is a closed expression for this integral; this former Macdonald
conjecture has been proved by Opdam, first for Weyl groups
\cite{someapplications} and later for finite reflection groups in general
\cite{besseletc}. For our needs it suffices to know that $c_k\neq 0$ if
$\Re{k}\geq 0$, which can be proved by more elementary means \cite[Corollary
4.17]{dunkltransform}.

The first results for the transform, and notably a version of the Plancherel
theorem, were obtained by Dunkl \cite{hankeltransform}. Later on, a more
systematic study was undertaken in \cite{dunkltransform}. The main properties and results
are as follows.

\begin{theorem}\label{thm:inventionesresults} Let $\Re{k}\geq 0$
and $\x\in\ga$. Then:
\begin{enumerate}
\item $\dk\tk f=M_{i\x^*}\dk f\,\,(f\in\schw(\ga))$. \item $\ek\tk
f=-M_{i\x^*}\ek f\,\,(f\in\schw(\ga))$. \item $\dk M_{i\x^*}f=-\tk\dk f
\,\,(f\in\schw(\ga))$. \item $\ek M_{i\x^*}f=\tk\ek f\,\,(f\in\schw(\ga))$.
\item $\dk$ is a linear homeomorphism of $\schw(\ga)$, with inverse $\ek$.
\item If $f\in\lone$ and $\dk f\in\lone$ then $\dk\ek f=\ek\dk f$\null$=f$ a.e.
\item If $k\geq 0$ then $\dk$ maps $L_1(\ga,w_k(x)\,dx)\cap
L_2(\ga,w_k(x)\,dx)$ into $L_2(\ga,w_k(x)\,dx)$, isometrically with respect to
the two-norm corresponding to $w_k(x)\,dx$, and extends uniquely from
$L_1(\ga,w_k(x)\,dx)\cap L_2(\ga,w_k(x)\,dx)$ to a unitary operator on
$L_2(\ga,w_k(x)\,dx)$.
\end{enumerate}
\end{theorem}

\begin{remark}
The definition of Dunkl operators can be generalized to complex
reflection groups \cite{dunklopdam}.
\end{remark}

\section{Formularium}\label{sec:formularium}

In this section we establish some useful formulas. Some of them are known from
the work of Dunkl, who based his proofs to a large extent on the existence of a
$k$-harmonic decomposition for polynomials \cite[Theorem
1.7]{reflectiongroups}. However, with the benefit of hindsight we can simplify
some of the original proofs considerably, by systematically exploiting the
following commutator relation \cite [Proposition 2.2]{diffdiff}:
\begin{equation}\label{eq:basisrelation}
\left[\mxi,\frac{\deltak}{2}\right]=-\tk\quad(\x\in\ga).
\end{equation}

To start, we first note that $\deltak:\poly\mapsto\poly$ is homogeneous of
degree $-2$, hence locally nilpotent. This enables us to define the linear
automorphisms $e^{\pm\frac{\deltak}{2}}$ of $\poly$. These automorphisms will
play an important part in what follows, together with the Gaussian
$\psi(x)=\exp(-\vert x\vert^2/2)$. It is known \cite{hankeltransform}
\cite{dunkltransform} that the Gaussian is an eigenfunction of both $\dk$ and
$\ek$ with eigenvalue $1$.

\begin{lemma}\label{lem:commutatorrelations}
Let $\x\in\ga$. Then for arbitrary $k$ we have
\begin{enumerate}
\item in ${\End}_\C(\poly)$:     $\left[\mxi\,,\,\expmindeltak\right]=\tk\circ\expmindeltak$;
\item in ${\Hom}_\C(\poly,\schw(\ga))$: $\tk\circ\mpsi\circ\expmindeltak=-\mpsi\circ\expmindeltak\circ\mxi$.
\end{enumerate}
\end{lemma}

\begin{proof}
The first part follows immediately from \eqref{eq:basisrelation} and the
obvious fact that $[\deltak,\tk]=0$. As to the second part, note that
$\tk\circ\mpsi=\mpsi\circ\tk-\mpsi\circ\mxi$, as a consequence of the
$G$-invariance of $\psi$. Hence
\begin{align*}
\tk\circ\mpsi\circ\expmindeltak&=\mpsi\circ\tk\circ\expmindeltak -
\mpsi\circ\mxi\circ\expmindeltak
\\
&=\mpsi\circ\tk\circ\expmindeltak-\mpsi\circ\left\{\left[\mxi\,,\,\expmindeltak
\right]+\expmindeltak\circ\mxi\right\}
\\
&=\mpsi\circ\tk\circ\expmindeltak-\mpsi\circ\left\{\tk\circ\expmindeltak +
\expmindeltak\circ\mxi\right\}
\\
&=-\mpsi\circ\expmindeltak\circ\mxi.
\end{align*}
\end{proof}
Repeated application of the second part of Lemma~\ref{lem:commutatorrelations}
yields the formula
\begin{equation}\label{eq:hompolyschwartzequation}
T_p(k) \left\{\left(\expmindeltak q\right)\psi\right\}=
(-1)^{\deg{p}}\left\{\expmindeltak (pq)  \right\}\psi,
\end{equation}
for homogeneous $p\in\poly$ and arbitrary $q\in\poly$. Taking $q=1$ this
implies, together with Theorem~\ref{thm:inventionesresults}, the following
result, which is equivalent to \cite[Proposition 2.1]{hankeltransform}.

\begin{corollary} \label{cor:transformofpolypsi} If $p\in\poly$ is homogeneous and $\Re k\geq 0$, then
\begin{enumerate}
\item
$
\dk(p\psi)=(-i)^{\deg{p}} \left(\expmindeltak p \right)\psi.
$
\item
$
\ek(p\psi)=i^{\,\deg{p}} \left(\expmindeltak p \right)\psi.
$
\end{enumerate}
\end{corollary}

The second part of Lemma~\ref{lem:commutatorrelations} also enables us to
reprove the symmetry of a bilinear form which was introduced by Dunkl
\cite{integralkernels}, as follows. For $p,q\in\poly$, put
$(p,q)_k=\left(T_p(k) q\right)(0)$. Although it is not obvious from the
definition, this bilinear form on $\poly$ is actually symmetric. This symmetry
follows from the following generalization by Dunkl [loc.cit.] of a result of
Macdonald:
\begin{equation}\label{eq:bilinearformequality}
(p,q)_k=\frac{1}{c_k}\inta \left(\expmindeltak p\right)\left(\expmindeltak q
\right)\psi w_k \,dx\quad(\Re k\geq 0, \,p,q\in\poly),
\end{equation}
in which the right hand side is obviously symmetric. In order to re-establish
\eqref{eq:bilinearformequality}, denote the right hand side by $[p,q]_k$. Now
\eqref{eq:antisymmetry} and Lemma~\ref{lem:commutatorrelations} imply that
\begin{align*}
[p,\tk q]_k&=\frac{1}{c_k}\inta \left(\expmindeltak p\right)\left(\tk\expmindeltak q
\right)\psi w_k\,dx
\\
&=\frac{1}{c_k}\inta -\tk\left\{\left(\expmindeltak p\right)\psi\right\}
\left(\expmindeltak q\right) w_k \,dx
\\
&=\frac{1}{c_k}\inta \left\{\expmindeltak (\mxi p)\right\}\left(\expmindeltak q
\right)\psi w_k \,dx\\
&=[\mxi p,q]_k.
\end{align*}
But the other form $(\,.\,,\,.\,)_k$ also has this property: $(\mxi
p,q)_k=(p,\tk q)_k$, as a direct consequence of its definition. Since it is
easy to see that $(1,q)_k=[1,q]_k$, an induction with respect to $\deg{p}$ then
proves that $(p,q)_k=[p\,,\,q]_k$, which is \eqref{eq:bilinearformequality}.

The following proposition will be used in the reduction of the proof of the
Paley--Wiener theorem~\ref{thm:pw} to the one-dimensional even case.

\begin{proposition}\label{prop:basisrelation} If $\Re k,\,\Re k^\prime
\geq 0$, then
\[
E_{k^\prime}D_k\left(p\psi\right)=\left(\expmindeltakp \expdeltak p\right)\psi
\quad(p\in\poly).
\]
\end{proposition}

\begin{proof}
We may assume that $p$ is homogeneous. Then Corollary~\ref{cor:transformofpolypsi} implies that
\begin{align*}
E_{k^\prime}\dk\left(p\psi\right)&=(-i)^{\deg p} E_{k^\prime} \left[ \left\{ \sum_{n=0}^\infty \frac{1}{n!} \left(-\frac{\deltak}{2}\right)^n p\right\} \psi\right]\\
&=(-i)^{\deg p} \left[ \sum_{n=0}^\infty  \frac{i^{\,\deg p -2n}}{n!} \expmindeltakp \left\{\left(-\frac{\deltak}{2}\right)^n p\right\} \right] \psi\\
&=\left(\expmindeltakp\expdeltak p\right)\psi.
\end{align*}
\end{proof}

To conclude this section, we re-establish a formula of Heckman
\cite{heckmanremark} which expresses the pivotal r\^ole of the $k$-Laplacian
$\D_k$:
\begin{equation}\label{eq:heckmanresult}
T_p(k)=\frac{1}{n!}\left(\ad \,\frac{\D_k}{2}\right)^n M_p,
\end{equation}
for $p\in\poly_n$. The proof in [loc.cit.] is based on representation theory
for $\textup{sl}(2)$, but it can also be seen directly, as follows. We may
assume that $M_p=\prod_{i=1}^{n}M_{\x_i^*}$. Then
\begin{equation*}
\left(\ad \,\frac{\D_k}{2}\right)^n M_p =\sum_{\sum j_i=n}\frac{n!}{j_1!\cdots
j_n!}\left\{ \left(\ad
\,\frac{\D_k}{2}\right)^{j_1}M_{\x_1^*}\right\}\circ\cdots\circ\left\{\left(\ad
\,\frac{\D_k}{2}\right)^{j_n}M_{\x_n^*}\right\}.
\end{equation*}
But $\left(\ad \,\D_k\right)^2M_{\x_i^*}=0$ as a consequence of
\eqref{eq:basisrelation} and the commutativity of the $\tk$, so the only
surviving term in the summation is the one with $j_1=\ldots=j_n=1$. Using
\eqref{eq:basisrelation} once more, this proves \eqref{eq:heckmanresult} .

\section{Paley--Wiener theorems}\label{sec:pw}

In this section, we conjecture a geometrical form of the Paley--Wiener theorem
for the Dunkl transform, and present several theorems to support it.

Establishing notation, for $S\subset\ga$ we let $\test(S)$ denote the smooth
compactly supported functions with support contained in $S$. If $S$ is compact
and non-empty, then we define the indicator $I_S:\ga\mapsto\R$ as
$I_S(x)=\max_{y\in S}(x,y)$ for $x\in\ga$. For such $S$, let ${\mathcal H}_S$
be the functions on $\gac$ of Paley--Wiener type corresponding to $S$, i.e.,
those entire functions with the property that for each integer $M\geq 0$ there
exists a constant $\g_M$ such that $\vert f(\lambda)\vert\leq \g_M
(1+\vert\lambda\vert)^{-M}\exp I_S(\Im{\lambda})$ for all $\lambda\in\gac$. We
then conjecture the following.

\begin{conjecture}[Paley--Wiener conjecture]\label{conj:pw}
Let $G$ be a finite reflection group. If $\Re k\geq 0$, and if $S$ is a
non-empty $G$-invariant compact convex subset of $\ga$, then $D_k$ is a linear
isomorphism between $\test(S)$ and ${\mathcal H}_S$.
\end{conjecture}

There is some evidence supporting this conjecture:
\begin{itemize}

\item The inversion theorem and \cite[Corollary 4.10]{dunkltransform} show that $D_k$ is
an injective map from $\test(S)$ into ${\mathcal H}_S$.

\item For $k=0$ the statement holds \cite[Theorem~7.3.1]{hormander}.

\item If $S$ is a ball centered at the origin, then a reduction to the
one-dimensional even case, where asymptotics can be used, enables us to
establish the statement as Theorem~\ref{thm:pw} below.

\item If the $k_\alpha$ are all strictly positive integers, then the statement
can be established, using shift operators, as Theorem~\ref{thm:convexpw} below.

\item If $G$ is a Weyl group, and if $S$ is the intersection of convex hulls of
orbits, then a limit transition from results of Opdam establishes the statement
as Theorem~\ref{thm:pwhecke} below.
\end{itemize}

The main obstacle for a possible proof of the conjecture along the usual lines,
using a contour shift, is the absence of adequate asymptotic results for the
Dunkl kernel. There are
some asymptotic results available \cite{roeslerdejeu}, but these fall far short
of what is needed. It is to be expected that better results could be obtained if
more were known about R\"osler's representing measures
\cite{roeslerpositivity}, but as yet these remain elusive. But even if much
stronger asymptotic results became available, the proof of
Theorem~\ref{thm:rankonecase} below seems to suggest that additional monodromy
arguments may then still be necessary.

\begin{remark}
In \cite{trimeche2}, a proof of Conjecture~\ref{conj:pw} if $k\geq 0$ and
$\sumposroots k_\alpha>0$ is presented. That proof, however, is not correct and to our knowledge Conjecture~4.1 is at the time of writing still open.
\end{remark}

\subsection{The case of arbitrary $G$ and $\Re k\geq 0$}\label{subsec:arbitrary}

Throughout this section, $\ball$ will denote the closed ball in $\ga$ with
radius $R$ and the origin as center. The space $\test(\ball)$ carries the usual
Fr\'echet topology of uniform convergence of all derivatives.

Our approach of the Paley--Wiener theorem for such sets---under the assumption
that $\Re k\geq 0$---consists of three steps, cf.\ the proof of Theorem~\ref{thm:pw}. First, we prove the result for
even functions in one dimension. Second, it is shown that this implies the
theorem for radial functions in arbitrary dimension. In the third step we
finally prove that the result for radial functions implies the theorem for
general smooth functions with support in a closed ball $\ball$.

\begin{remark}
The special r\^ole of radial functions in the theory has been noted by several
authors \cite{thesis} \cite{xu} \cite{roeslervoitmarkov} \cite{trimeche1},
e.g., the Dunkl transform of a radial function is again radial
\cite{roeslervoitmarkov}. In Paley--Wiener theory this phenomenon is encountered once more, when reducing the radial case in arbitrary dimension to the
even case in one dimension. The one-dimensional
even case can be handled by either Weyl fractional integral operators
\cite{trimeche1}, or as in our approach below, which uses classical results for
the asymptotics of Bessel functions combined with a contour shift. The final
step in our approach, from the radial to the general case, and for which a reference to
\cite{thesis} was given in \cite{trimeche1}, is in a sense the deepest, since
it ultimately rests on the general theory of representations of compact groups
in Fr\'echet spaces. Recently, an alternative reduction to the one-dimensional case has been introduced which circumvents the inference of the general case from the radial one \cite{thangaveluxu}. This reduction is based on results for $k$-harmonic polynomials and for functions supported in $B_R$ it gives a proof of the Paley--Wiener theorem which is independent of the results in the present paper.

\end{remark}

Starting with the proof, the one-dimensional even case is settled in the
following theorem, where, as usual, invariance is denoted by superscripts.

\begin{theorem}\label{thm:rankonecase} Let $\ga$ be one-dimensional and $\Re k\geq 0$. If $f\in{\mathcal H}_{\ball}^{\Z_2}$, then $E_k f\in\test(\ball)^{\Z_2}$.
\end{theorem}

\begin{proof}
Let $f$ be as in the statement. Then obviously $E_k f$ is an even rapidly
decreasing smooth function. We have to show that
\begin{equation*}
\int_{-\infty}^{\infty} f(\lambda)\textup{Exp}_{\Z_2} (i\lambda,k,x)
\vert\lambda\vert^{2k}\,d\lambda=0
\end{equation*}
if $x>R$. Fix such $x$. From \cite{hankeltransform} we have
\[\textup{Exp}_{\Z_2} (i\lambda,k,x)=\Gamma (k+\onehalf) \left(\frac{\lambda
x}{2}\right)^{\onehalf-k} \left\{J_{k-\onehalf}(\lambda
x)+iJ_{k+\onehalf}(\lambda x)\right\}.\] Using the invariance of $f$ we see
that we are left with showing that
\begin{equation}\label{eq:rankoneintegral}
\int_{-\infty}^\infty f(\lambda) \left({\lambda x}\right)^{\onehalf-k}
J_{k-\onehalf}(\lambda x)\vert \lambda\vert ^{2k}\,d\lambda=0.
 \end{equation}
This expression makes it obvious that there are two obstructions for a direct
application of the classical argument of shifting the contour to infinity.
First, the Bessel function has exponential growth in both the positive and
negative imaginary directions, and second the weight function
$\vert\lambda\vert^{2k}$ is in general not the restriction of a holomorphic
function on the upper or lower half plane. So we proceed indirectly.

Recall the definition of the Bessel functions of the third kind:
\begin{align*}
H_\n^{(1)}&=\frac{J_{-\n} - e^{-\n\p i}J_\n}{i\sin \n\p},\\
H_\n^{(2)}&=\frac{J_{-\n} - e^{\n\p i}J_\n}{-i\sin \n\p},
\end{align*}
so that
\[
J_\n=\frac{H_\n^{(1)}+H_\n^{(2)}}{2}.
\]
If $\n$ is an integer then a limit has to be taken. In our case, this occurs if
$k$ is half-integer, but by continuity in $k$ we may and will assume that this
is not the case.

For our purposes, the important property of these functions is the
asymptotic behaviour \cite[9.2.7]{abramowitz}:
\begin{equation}\label{eq:hankelasymptotics}
H_\n^{(1)}(z)=\sqrt{\frac{2}{\p z}} e^{i\left(z-\frac{\n\p}{2}-\frac{\p}{4}\right)} \left(1+O\left(\frac{1}{z}\right)\right),
\end{equation}
valid if $-\p<\arg \,z<2\p$ (which is to be interpreted in the sense of
analytic continuation). Note that the range of validity of this asymptotic
development contains the entire upper half plane and that (in contrast to the
ordinary Bessel function) this Hankel function has exponential decrease in the positive
imaginary direction.

Define $\phi^{(1)},\phi^{(2)}:(0,\infty)\mapsto\C$ by
$$
\phi^{(1)}(\lambda)=(\lambda x)^{\onehalf-k}H_{k-\onehalf}^{(1)}(\lambda
x)\lambda^{2k}
$$
and
$$
\phi^{(2)}(\lambda)=(\lambda x)^{\onehalf-k}H_{k-\onehalf}^{(2)}(\lambda
x)\lambda^{2k}.
$$
Let $\phi^{(1)}_{\tu{c}}$ denote the analytic continuation of $\phi^{(1)}$ from
$(0,\infty)$ to $\C\setminus\{it\mid t\leq 0\}$. If one recalls that
$J_\n(z)=z^\n \tilde{J}_\n (z)$ with $\tilde{J}_\n$ entire and
$\Z_2$-invariant, one notes that $\phi^{(1)}_{\tu{c}}(\lambda)$ remains bounded as
$\lambda\rightarrow 0$ in $\C\setminus\{it\mid t\leq 0\}$, since $\Re k\geq 0$.
A small computation will also make it clear that $\phi^{(1)}_{\tu{c}}
(\lambda)=\phi^{(2)}(-\lambda)$ for $\lambda<0$.

Using the invariance of $f$ and the weight function, we then compute as
follows.
\begin{align*}
\int_{-\infty}^\infty f(\lambda) \left({\lambda x}\right)^{\onehalf-k}
&J_{k-\onehalf}(\lambda x)\vert \lambda\vert ^{2k}\,d\lambda=\\&=2\int_0^\infty
f(\lambda) \left({\lambda x}\right)^{\onehalf-k} J_{k-\onehalf}(\lambda
x)\lambda ^{2k}\,d\lambda\\
&=\int_0^\infty f(\lambda) \left({\lambda
x}\right)^{\onehalf-k}\left(H_{k-\onehalf}^{(1)}(\lambda
x)+H_{k-\onehalf}^{(2)}(\lambda x)\right)\lambda^{2k}\,d\lambda\\
&=\int_0^\infty f(\lambda)\left(\phi^{(1)}(\lambda)+
\phi^{(2)}(\lambda)\right)\,d\lambda\\
&= \lim_{\epsilon\downarrow 0}\left\{\int_\epsilon^\infty f(\lambda)
\phi^{(1)}(\lambda)\,d\lambda + \int_\epsilon^\infty f(\lambda)
\phi^{(2)}(\lambda)\,d\lambda\right\}\\
&=\lim_{\epsilon\downarrow 0}\left\{\int_\epsilon^\infty f(\lambda)
\phi^{(1)}(\lambda)\,d\lambda + \int_{-\infty}^{-\epsilon} f(-\lambda)
\phi^{(2)}(-\lambda)\,d\lambda\right\}\\
&= \lim_{\epsilon\downarrow 0}\left\{\int_\epsilon^\infty f(\lambda)
\phi^{(1)}_{\tu{c}}(\lambda)\,d\lambda + \int_{-\infty}^{-\epsilon} f(\lambda)
\phi^{(1)}_{\tu{c}}(\lambda)\,d\lambda\right\}.
\end{align*}
Now $f\phi^{(1)}_{\tu{c}}$ is holomorphic on $\C\setminus\{it\mid t\leq 0\}$, and,
since $x>R$, it has exponential decrease in the positive imaginary direction,
as a consequence of \eqref{eq:hankelasymptotics}. The classical argument
therefore establishes \eqref{eq:rankoneintegral}, with a minor modification
involving a semi-circle of radius $\e$ around $0$ in the upper half plane and
using that $f\phi^{1}_{\tu{c}}$ is bounded around $0$.
\end{proof}

Next, we proceed with the reduction of the general radial case to the
one-dimensional even case. To this end, fix $x_0\neq 0$ in $\ga$, and define
the map $\Resxo:\schwinv\mapsto\schw(\R)^{\Z_2}$ by restricting to the line
passing through $x_0$:
\[
(\Resxo f)(s)=f\left(s\frac{x_0}{\vert
x_0\vert}\right)\quad\left(s\in\R,\,f\in\schwinv\right).
\]

The following proposition, which implies the Paley-Wiener theorem for radial
functions (as will become apparent in the proof of Theorem~\ref{thm:pw}),
involves the Dunkl transform for general $\ga$ and $\R$ both at the same time.
We therefore add a subscript $\Z_2$ in the latter case for clarity. Let
$\schwinv$ denote the $\on$-invariants in $\schw (\ga)$.

\begin{proposition}\label{prop:rankonereduction}
Suppose $\Re k\geq 0$. Let $f\in\schwinv$ and put $\g=\sumposroots
k_\alpha$. Then
\[
\Resxo E_k D_0 f = E_{\g,\Z_2} D_{0,\Z_2}\Resxo f
\]
for all nonzero $x_0\in\ga$.
\end{proposition}

\begin{proof}
Let $\schwdense=\{p\psi\mid p\in\poly\}$, with $\psi$ denoting the Gaussian as
in Section~\ref{sec:formularium}. It is known that $\schwdense$ is dense in
$\schw(\ga)$, see \cite[p.~263]{distributions} for this particular result, or
\cite{densesubspaces} for a general framework for this type of problems. Note that
the canonical projection from $\schw(\ga)$ onto $\schwinv$ is continuous as a
consequence of the closed graph theorem, implying that $\schwinvdense$ is dense
in $\schwinv$. It is therefore, by linearity and continuity, sufficient to prove the
proposition for a function $f$ of the form $\vert x\vert^{2q} \psi$, where $q$
is a non-negative integer. This can be done using
Proposition~\ref{prop:basisrelation} and an identity for Laguerre polynomials.
Recall the definition:
\[
L_n^{(\alpha)}(x)=\sum_{m=0}^n \frac{(-1)^m}{m!} {n+\alpha\choose n-m} x^m.
\]
Then the following identity holds \cite[22.12.6]{abramowitz}:
\begin{equation}\label{eq:laguerreidentity}
\sum_{m=0}^n L_m^{(\alpha)}(x) L_{n-m}^{(\beta)}(y)=L_n^{(\alpha+\beta+1)}(x+y).
\end{equation}
It is known \cite[Proposition 3.9]{integralkernels} that $\expdeltazero \vert
x\vert^{2q}=2^q q! L_q^{(N/2\,-1)} \left(-\vert x\vert^2/2 \right)$ and
$\expmindeltak \vert x\vert^{2q}=(-2)^q q! L_q^{(N/2\,+\g-1)}\left(\vert
x\vert^2/2\right)$. Using Proposition~\ref{prop:basisrelation} we therefore
find
\begin{align*}
E_k D_0 \left(\vert x\vert^{2q}\psi\right)&=\left(\expmindeltak
\expdeltazero \vert x\vert^{2q}\right)\psi\\
&=2^q q!\sum_{m=0}^q (-1)^m {N/2\,+q-1 \choose q-m} L_m^{(N/2\,+\g-1)} \left(\vert x\vert^2/2\right)\psi\\
&=(-2)^q q!\sum_{m=0}^q {-N/2 -m \choose q-m} L_m^{(N/2\,+\g-1)} \left(\vert x\vert^2/2\right)\psi.
\end{align*}
The special case $n=q$, $y=0$, $\alpha=N/2 + \g -1$ and $\beta=-N/2 -q$ of \eqref{eq:laguerreidentity} then shows that
\[
E_k D_0 \left(\vert x\vert^{2q}\psi\right)=(-2)^q q! L_q^{(\g-q)}\left(\vert x\vert^2/2\right)\psi.
\]
Curiously enough, the dimension $N$ has dropped out, and the proposition
follows immediately from this observation.
\end{proof}

It follows easily from this result that the Dunkl transform of a radial
function is again radial, retrieving this result from \cite{roeslervoitmarkov}.

Let $\testbdense=\Span\{T_p(0)f\mid  p \in\poly,\,\,f\in\testbinv\}$. The following density result is crucial in the step from the radial to the general case.

\begin{proposition}\label{prop:densityproposition}
$\testbdense$ is dense in $\testb$.
\end{proposition}

This proposition implies a special case of \cite[Cor.~7.8, p.~310]{helref2},
but the latter result does not seem to imply the proposition.
The proof of Proposition~\ref{prop:densityproposition} is based on representation theory for compact groups in general and $\on$
in particular; we recall a few facts to start with.

Let $\p:K\mapsto\Aut(E)$ be a strongly continuous representation of the compact
group $K$ in the automorphism group of a Fr\'echet space $E$. A vector $e\in E$
is \emph{$K$-finite} if the orbit of $e$ under $K$ spans a finite dimensional
subspace of $E$. If $\delta\in\khat$ (the unitary dual of $K$), then a vector
$e\in E$ is \emph{$K$-finite of type $\delta$} if the span of the orbit of $e$
decomposes into finitely many copies of $\delta$. For $\delta\in\khat$, let
$E_\delta$ denote the subspace of $K$-finite vectors of type $\delta$. By
\cite[Lemma~1.9, p.~396]{helref} the subspace
$\bigoplus_{\delta\in\khat}E_\delta$ of $K$-finite vectors is dense in $E$.

The natural action of $\on$ on $\testb$ is such a strongly continuous
representation. Hence Proposition~\ref{prop:densityproposition} will be implied
by the density of the $K$-finite vectors, as mentioned above, once we have
proved the following result.

\begin{proposition}\label{prop:finitevectors}
$\testbdense$ is the space of $\on$-finite vectors in $\testb$.
\end{proposition}

Obviously, any element of $\testbdense$ is $\on$-finite, but the reverse
appears to be harder. The proof of Proposition~\ref{prop:finitevectors} below
uses some Fourier analysis; for clarity of notation, we will temporarily denote
the ordinary Fourier transform by $\fourier$ rather than $D_0$. We need the
following classical result in the representation theory of $\on$.

\begin{theorem}\label{thm:onrepresentationresults}
Let $\harm_n$ be the harmonic polynomials in $\poly_n$.
Any $p\in\poly_n$ has a unique decomposition of the form
\[
p=\sum_{l=0}^{[\frac{n}{2}]}\vert x\vert^{2l}p_l
\]
with $p_l\in \harm_{n-2l}$. The map $\p_{n,l}:\poly_n\mapsto\harm_{n-2l}$ that
sends $p$ to $p_l$ is $\on$-equivariant. Furthermore, the representations of
$\on$ on the spaces $\harm_n\,\,(n=0,1,\dots)$ are irreducible and mutually
inequivalent.
\end{theorem}

As a further preparation for the proof of Proposition~\ref{prop:finitevectors}
we need the following result.

\begin{lemma}\label{lem:onfourierversion}Let $\delta\in\ondual$ and $f\in\testb$. Suppose that the
orbit of $f$ under $\on$ spans a copy $V_\delta$ of $\delta$ in $\testb$. Then
there exist a harmonic homogeneous polynomial $p$ of type $\delta$ and scalars
$\lambda_m\,\,(m=0,1,\dots)$, such that $\fourier(f)(z)=
p(z)\sum_{m=0}^\infty\lambda_m(z,z)^m$ $(z\in\C^N)$, where the series converges
uniformly on compact subsets of $\C^N$.
\end{lemma}

\begin{proof}
By the Paley--Wiener theorem the Fourier transform $\fourier(f)$ is an entire
function on $\C^N$. $\fourier$ commutes with the $\on$-action, so there is a
copy of $\delta$ in the space of functions of Paley--Wiener type. Since the
action of $\on$ on entire functions is homogeneous (in the sense that the
action is on each homogeneous component of the power series development around $0$
separately), we conclude that $\delta$ also occurs in $\poly$. By
Theorem~\ref{thm:onrepresentationresults} there exists a unique copy, $\harm_l$
say, of $\delta$ as the harmonics of some homogeneous degree.

For $n\geq 0$, let $\p_{n,l}$ denote the map $\poly_n\mapsto \harm_l$ defined
by the decomposition in Theorem~\ref{thm:onrepresentationresults}. Also, for
$n\geq 0$ and $g\in V_\delta$, let $\fourier_n(g)$ denote the component of
homogeneous degree $n$ in the power series development of $\fourier(g)$ around
$0$. Consider the $\on$-equivariant maps $\p_{n,l}\circ
\fourier_n:V_\delta\mapsto \harm_l\,\,(n\geq 0)$. Since $\fourier(f)\neq 0$,
there exists $n_0$ such that $\p_{n_0,l}\circ\fourier_{n_0}$ is non-zero, hence
an isomorphism. But then for all $n$ the map
$(\p_{n,l}\circ\fourier_n)\circ(\p_{n_0,l}\circ\fourier_{n_0})^{-1}:
\harm_l\mapsto \harm_l$ must be equal to a multiple $\lambda_n$ of the
identity. Hence $f(z)=\sum_{m=0}^\infty \lambda_m (z,z)^m p(z)\,\,(z\in\C^N)$
for some nonzero $p\in\harm_l$. Since there exist $z\in\C^N$ with $p(z)\neq 0$
and $\vert z\vert$ arbitrarily large, the factor $p(z)$ can be brought outside
the summation.
\end{proof}

\begin{proof}[Proof of Proposition~\ref{prop:finitevectors}]
As remarked before, $\testbdense$ is certainly contained in the space of
$\on$-finite vectors in $\testb$. Conversely, let $f\in\testb$ be an
$\on$-finite vector. We may assume that the orbit of $f$ spans a copy of
$\delta$ in $\testb$ for some $\delta\in\ondual$. Applying
Lemma~\ref{lem:onfourierversion} one obtains a factorization $\fourier(f)=ph$
of $\fourier(f)$ in a polynomial $p$ and an entire $\on$-invariant function
$h$. We proceed to show that $h$ is in fact in $\mathcal H_{\ball}$. By
\cite[Lemma 8.3]{ruref} there exists a constant $A>0$ such that
\begin{equation}\label{eq:rudinintegral}
\vert h(z)\vert\leq A\int_{T^N}\vert\fourier(f)(z+w)\vert\,d\sigma_N(w)
\quad(z\in\C^N),
\end{equation}
where $T^N=\{(e^{i\theta_1} ,\dots,e^{i\theta_N} ) \mid
\theta_j\in\R\,\,(j=1,\dots,N)\}$ is the $N$-torus and $\sigma_N$ is the
normalized Haar measure on $T^N$. There exist constants
$\g_M\,\,(M=0,1,2,\dots)$ such that
\begin{equation}\label{eq:pwinequalities}
\vert\fourier(f)(z)\vert\leq\g_M(1+\vert z\vert)^{-M} e^{R\vert\Im z\vert}
\,\,(z\in\C^N).
\end{equation}
Now note that
\[
\frac{1+\vert z\vert}{1+\vert z+w\vert}\leq 1+ \frac{\vert w\vert}{1+\vert z+
w\vert}\leq 1+\sqrt {N}\quad(z\in\C^N,w\in T^N).
\]
Since in addition $\vert \Im (z+w)\vert\leq \vert \Im z\vert+\sqrt {N}$ for $z\in\C^N$ and $w\in T^N$, we
conclude from \eqref{eq:rudinintegral} that $h$ also satisfies estimates as in
\eqref{eq:pwinequalities}, i.e., $h\in\mathcal H_{\ball}$. Hence
$\fourier^{-1}(h)$ is in $\testb$, and in fact $\fourier^{-1}h\in\testbinv$
since $h$ is $\on$-invariant. This implies that $f\in\testbdense$, as required.
\end{proof}

This concludes the proof of Proposition~\ref{prop:densityproposition}.

We can now prove the Paley--Wiener theorem of this section by putting the
pieces together.

\begin{theorem}[Paley--Wiener theorem, first version]\label{thm:pw} Let $G$ be a finite reflection group  and suppose that $\Re k\geq 0$. Then the Dunkl transform $\dk$ is a linear isomorphism between $\testb$ and $\pwtype$, for all $R\geq 0$.
\end{theorem}

\begin{proof}
In view of the inversion theorem, all that needs to be done is to show that
$\ek g\in\testb$ for all $g\in\pwtype$. To this end, start by noting that $\ek
D_0 f\in\testb$ if $f\in\testbinv$, as an immediate consequence of
Proposition~\ref{prop:rankonereduction} and Theorem~\ref{thm:rankonecase}. But
then we also have $\ek D_0 f\in\testb$ for $f\in\testbdense$, since $\ek\circ
D_0\circ\partial_\x=\tk\circ\ek\circ D_0$. From
Proposition~\ref{prop:densityproposition} and continuity we can then conclude
that $\ek D_0 f\in\testb$ for $f\in\testb$. By the Paley--Wiener theorem for
$D_0$, we see that $\ek g \in\testb$ for all $g\in\pwtype$, as required.
\end{proof}

\subsection{The case of arbitrary $G$ and strictly positive integral $k$}

If all $k_\alpha$ are strictly positive integers, Conjecture~\ref{conj:pw} can
be established from the corresponding theorem for the Fourier transform, by
showing that $\Exp$ is then in fact an elementary function. This can be seen by
starting from Heckman's work on shift operators \cite{heckmanremark}.

We introduce some notation. Let $R=\bigcup_{i=1}^m S_i$ be the decomposition of
$R$ into $G$-orbits, and put $S_{i,+}=S_i\cap\rplus$. Define
$\p_i=\prod_{\alpha\in S_{i,+}}\alpha^*$, and let Res denote restriction to the
$G$-invariant functions. Then we have the following results \cite{heckmanremark}, valid for
$p\in\poly^G$:

\begin{align}
\Res \left(M_{\p_i}^{-1}T_{\p_i}(k)\right) \Res \left(T_p(k)\right)&=\Res \left(T_p(k+1_{S_i})\right)
\Res \left(M_{\p_i}^{-1}T_{\p_i}(k)\right),
\label{eq:firstshiftrelation}\\
\Res \left(T_{\p_i}(k)M_{\pi_i}\right)
\Res \left(T_p(k+1_{S_i})\right) &= \Res \left(T_p(k)\right)
\Res \left(T_{\p_i}(k)M_{\pi_i}\right),\label{eq:secondshiftrelation}
\end{align}
where $(k+1_{S_i})_\alpha=k_\alpha+1$ if $\alpha\in S_i$ and  $(k+1_{S_i})_\alpha=k_\alpha$ if $\alpha\notin S_i$.

The generalized Bessel kernel $J_G$ is defined for $k\in\kreg$ by
\begin{equation}\label{eq:besselkerneldef}
J_G(\lambda,k,x)=\frac{1}{\vert G\vert}\sum_{g\in G} \Exp (\lambda,k,gx)
\end{equation}
and is (up to multiples) the unique $G$-invariant solution of the Bessel system
\[
T_p(k)f=p(\lambda)f\quad(p\in\poly^G).
\]
It follows immediately from \eqref{eq:firstshiftrelation},
\eqref{eq:secondshiftrelation} and the uniqueness of invariant solutions of the
Bessel system that there exist complex constants $c_{1,2}(\lambda,k)$ such that
\[
M_{\p_i}^{-1}T_{\p_i}(k) J_G(\lambda,k,\,.\,)= c_1(\lambda,k)J_G(\lambda,k+1_{S_i},\,.\,)
\]
and
\[
T_{\p_i}(k)M_{\pi_i} J_G(\lambda,k+1_{S_i},\,.\,)= c_2(\lambda,k)J_G(\lambda,k,\,.\,).
\]
Evaluation at zero yields $c_2(\lambda,k)=(\p_i,\p_i)_k$, where $(\,.\,,\,.\,)_k$ is the bilinear form from Section~\ref{sec:formularium}.
\begin{align*}
c_1(\lambda,k)(\p_i,\p_i)_k J_G(\lambda,k,\,.\,)&=c_1(\lambda,k)
T_{\p_i}(k)M_{\pi_i} J_G(\lambda,k+1_{S_i},\,.\,)
\\
&= T_{\p_i}(k)M_{\pi_i}  M_{\p_i}^{-1}T_{\p_i}(k) J_G(\lambda,k,\,.\,)
\\
&=  \p_i^2(\lambda)J_G(\lambda,k,\,.\,).
\end{align*}
Hence
\begin{equation}\label{eq:basicbesselrelation}
\p_i^2(\lambda) J_G(\lambda,k+1_{S_i},\,.\,) = (\p_i,\p_i)_k M_{\p_i}^{-1}T_{\p_i}(k) J_G(\lambda,k,\,.\,).
\end{equation}
Suppose then that the (fixed) value of $k$ on $S_i$ is $k_i$, a strictly
positive integer. Repeated application of \eqref{eq:basicbesselrelation} yields
\begin{align}
\label{eq:complicatedformula}
w_k(\lambda)J_G(\lambda,k,\,.\,) =&{\rm~} \prod_{i=1}^m \prod_{n=1}^{k_i} (\p_i,\p_i)_{(n-1)1_{S_i} + \sum_{j<i}k_j 1_{S_j}}\cdot\\
& \cdot  \prod_{i=1}^m \prod_{n=1}^{k_i} M_{\p_i}^{-1} T_{\p_i}
\left((n-1)1_{S_i} + \sum_{j<i}k_j 1_{S_j}\right)\,\, \frac{1}{\vert G\vert}
\sum_{g\in G} e^{(g\lambda,\,.\,)}.\notag
\end{align}
We conclude from \eqref{eq:complicatedformula} that there exist functions
$p_g(\lambda,x)$ on $\gac\times\gacreg$ (we suppress the dependence on our
fixed $k$) which are polynomials in $\lambda$ with rational functions in $x$ as
coefficients (with poles with respect to $x$ along the singular hyperplanes, if
any) such that
\[
w_k(\lambda) J_G(\lambda,k,x)=\sum_{g\in G} p_g(\lambda,x) e^{(g\lambda,x)}\quad (\lambda\in\gac,\,x\in\ga^{\reg}),
\]
where as usual $\gareg=\{x\in\ga\mid (\alpha,x)\neq 0\,\,\forall\alpha\in R\}$
and $\gacreg$ is defined similarly. Now the first operator acting in
\eqref{eq:complicatedformula} is $T_{\p_1}(0)$, showing that the $p_g$ are
divisible by $\p_1(\lambda)$. Since the numbering of the $S_i$ is arbitrary we
see that the $p_g$ are in fact divisible by $\p(\lambda)=\prod_{i=1}^m
\p_i(\lambda)$. Thus
\begin{equation}\label{eq:besselelementary}
w_k(\lambda) J_G(\lambda,k,x)=\pi(\lambda) \sum_{g\in G} p_g^\prime(\lambda,x) e^{(g\lambda,x)}\quad (\lambda\in\gac,\,x\in\gareg),
\end{equation}
where the $p_g^\prime$ are of the same form as the $p_g$. It is here that we
have used that the $k_\alpha$ are \emph{strictly} positive integers, in order
to obtain the factor $\pi(\lambda)$ in the right hand side of
\eqref{eq:besselelementary} which is needed below.

There is a similar expression for $\Exp$, which we obtain using a result of
Opdam \cite{besseletc} describing how $\Exp$ may be obtained from $J_G$.
According to [loc.cit., 6.4-6.6] there exists a rational function
$Q(\lambda,x)$ on $\gac\times\gac$ such that:
\begin{enumerate}
\item[$\bullet$] $(\lambda,x)\mapsto \p(\lambda)Q(\lambda,x)$ is a polynomial
on $\gac\times\gac$; \item[$\bullet$]
$T_{Q(\lambda,\,.\,)}(k)J_G(\lambda,k,\,.\,)=\Exp (\lambda,k,\,.\,)
\,\,(\lambda\in\gacreg,\,x\in\gac)$.
\end{enumerate}
Using this and \eqref{eq:besselelementary} we see that, for $\lambda\in\gacreg$ and $x\in\gareg$,
\begin{align}\label{eq:expelementary}
w_k(\lambda)\Exp(\lambda,k,x)&=T_{Q(\lambda,\,.\,)} w_k(\lambda) J_G(\lambda,k,\,.\,)\\
&=T_{Q(\lambda,\,.\,)}\pi(\lambda)\sum_{g\in G} p_g^{\prime}(\lambda,x) e^{(g\lambda,x)}\notag\\
&=T_{\pi(\lambda)Q(\lambda,\,.\,)}\sum_{g\in G} p_g^{\prime}(\lambda,x) e^{(g\lambda,x)}\notag\\
&=\sum_{g\in G} p_g^{\prime\prime}(\lambda,x) e^{(g\lambda,x)},\notag
\end{align}
where again the $p_g^{\prime\prime}$ are of the same form as the $p_g$. But then \eqref{eq:expelementary} is by continuity
in fact valid for $\lambda\in\gac$ (rather than $\lambda\in\gacreg$) and
$x\in\gareg$. This observation allows us to establish Conjecture~\ref{conj:pw}
for strictly positive integral multiplicities.

\begin{theorem}[Paley--Wiener theorem, second version]\label{thm:convexpw}
Let $G$ be a finite reflection group and suppose that the $k_\alpha$ are all
strictly positive integers. If $S$ is a non-empty $G$-invariant
compact convex subset of $\ga$, then $D_k$ is a linear isomorphism between
$\test(S)$ and ${\mathcal H}_S$.
\end{theorem}

\begin{proof}
Suppose $f\in{\mathcal H}_S$, $x\in\ga$, $x\notin S$. In view of the inversion
theorem, the only thing that remains to be proved is:
\begin{equation}\label{eq:integraliszero}
\inta f(\lambda)\Exp (i\lambda,k,x)w_k(\lambda)\,d\lambda=0.
\end{equation}
If $x\in\gareg$, then \eqref{eq:expelementary} allows us to write this integral as
$$
\sum_{g\in G}\inta f(\lambda) p_g^{\prime\prime}(i \lambda,x) e^{i(\lambda,g^{-1} x)} \,d\lambda,
$$
where the $p_g^{\prime\prime}$ are polynomials in $\lambda$ with rational
functions in $x$ as coefficients (with poles with respect to $x$ along the
singular hyperplanes, if any). Hence the geometrical form of the Paley--Wiener
theorem for the Fourier transform \cite[Theorem~7.3.1]{hormander} shows that
each of the summands is zero. Then \eqref{eq:integraliszero} holds for all
$x\notin S$ by continuity.
\end{proof}

\subsection{The case where $G$ is a Weyl group and $\Re k\geq 0$}\label{subsec:gradhecke}

In the case where $G$ is a Weyl group and $\Re k\geq 0$, a Paley--Wiener
theorem for intersections of convex hulls of orbits can be obtained by a limit
transition from the results of Opdam on Cherednik operators for $k\geq 0$
\cite{gradhecke}. If every non-empty $G$-invariant compact convex set were the
intersection of the convex hulls of orbits, then this would establish
Conjecture~\ref{conj:pw} for Weyl groups. However, this intersection property
does not hold, as the two-dimensional example of the closed unit ball and
$A_1\times A_1$ shows. Thus, Conjecture~\ref{conj:pw} remains open even
for Weyl groups, and this example also shows that in this context neither one of the
Theorems~\ref{thm:pw} above and \ref{thm:pwhecke} below implies the other.

We start by collecting the relevant results from \cite{gradhecke}.

Suppose that $R$ is an integral root system, and let $k:R\mapsto\C$ be a multiplicity function. We will assume that $R$ is reduced, which is sufficient for our purposes. Choose a system $\rplus$ of positive roots, and define the Cherednik operator $D_\xi$, for $\xi\in\gac$, as
\[
D_\xi(k)=\partial_\xi + \sumposroots k_\alpha (\alpha,\xi) M_{\textstyle{(1-e^{-\alpha^*})^{-1}}} (1-\refl) - (\rho(k),\xi),
\]
where $\rho(k)=\onehalf\sum_{\alpha\in\rplus} k_\alpha\alpha$. The Cherednik
operators are not equivariant, and their definition is dependent on the choice
of $\rplus$. For fixed multiplicity they form a commutative family, and together with $G$ they
then generate an algebra of operators which is an isomorphic copy of the graded
Hecke algebra, corresponding to the choice of $\rplus$ and $k$, as it is
defined by Lusztig \cite{lusztig}. We refer to \cite{cherednik} for this
isomorphism, or to \cite{gradhecke}. Aside, we remark that the $D_\xi(k)$
certainly depend on the length of the roots, in contrast to the Dunkl
operators, but it will become apparent that this is innocent for our purposes.

The Dunkl kernel $\Exp$ can be obtained from the eigenfunctions of the Cherednik operators by taking a suitable limit. To see this, we recall a weakened version of some of Opdam's results \cite[c.f.\ Theorem~3.15 and Proposition~6.1]{gradhecke}.

\begin{theorem}
For all $k\geq 0$ there exists an open neighborhood $U$ of $0\in\ga$ and a
holomorphic function $C_G(\,.\,,k,\,.\,):\gac\times(\ga + iU)\mapsto\C$ with
the following properties:
\begin{enumerate}
\item $C_G(\lambda,k,0)=1$ for all $\lambda\in\gac$; \item
$D_\xi(k)C_G(\lambda,k,\,.\,)=(\lambda,\xi)C_G(\lambda,k,\,.\,)$ for all
$\lambda,\,\xi\in\gac$; \item if $\lambda\in\gac$, and $z=x+iy$ with $x\in\ga$
and $y\in U$, then
\begin{equation}\label{eq:cherednikestimate}
\vert C_G(\lambda,k,z)\leq\sqrt{\vert G\vert} e^{-\min_{g\in
G}(g\,\Im\lambda,y)+\max_{g\in G}(g\,\rho(k),y)+\max_{g\in
G}(g\,\Re\lambda,x)}.
\end{equation}
\end{enumerate}

\end{theorem}
The estimate in \eqref{eq:cherednikestimate} enables us to prove the following limit transition.

\begin{theorem}\label{thm:limittransition}
Let $k\geq 0$ be fixed, and suppose that $S\subset\gac\times\gac$ is a non-empty compact subset. Then, as $\epsilon\rightarrow 0$ through the complex numbers, we have
\begin{equation}\label{eq:limittransition}
\lim_{\epsilon\rightarrow 0} C_G(\epsilon^{-1}\lambda,k,\epsilon z)=\Exp(\lambda,k,z),
\end{equation}
uniformly for $(\lambda,z)\in S$. Here it is understood that $\vert \epsilon\vert$ is taken sufficiently small, so that $C_G(\epsilon^{-1}\lambda,k,\epsilon z)$ is defined for all $(\lambda,z)\in S$.
\end{theorem}

\begin{proof}
Suppose, to the contrary, that there exist $\eta>0$, a non-empty compact set $S\subset\gac\times\gac$, and a complex sequence $\epsilon_n\rightarrow 0$ such that, for all $n$,
\begin{equation}\label{eq:nonconvergence}
\sup_{(\lambda,z)\in S}\vert C_G(\epsilon_n^{-1}\lambda,k,\epsilon_n z)-\Exp (\lambda,k,z)\vert\geq\eta.
\end{equation}
Take an open neighborhood of $S$ of the form $V_1\times V_2$, with
$V_1,V_2\subset\gac$ open balls centered at the origin, and define
$\phi_n:V_1\times V_2\mapsto\C$ by
$\phi_n(\lambda,z)=C_G(\epsilon_n^{-1}\lambda,k,\epsilon_n z)$, discarding a
finite number of $\epsilon_n$ if necessary. From \eqref{eq:cherednikestimate}
one sees that the sequence $\{\phi_n\}_{n=1}^\infty$ is uniformly bounded on
$V_1\times V_2$. By Montel's theorem we may therefore assume, passing to a
subsequence if necessary, that the $\phi_n$ converge uniformly on compact
subsets of $V_1\times V_2$ to a holomorphic function $\phi_\infty$ on
$V_1\times V_2$.

Now a small computation shows, for fixed $\lambda\in V_1$ and for arbitrary $\xi\in\gac$, that as functions on $V_2$ we have
\[
\left\{\partial_\xi + \epsilon_n \sumposroots k_\alpha \frac{(\alpha,\xi)}{1-e^{-\epsilon_n(\alpha,\,.\,)}} (1-\refl)\right\}\phi_n(\lambda,\,.\,)=\left\{(\xi,\lambda)+\epsilon_n(\rho(k),\xi)\right\}\phi_n(\lambda_,\,.\,).
\]
Passing to the limit one obtains
\[
\left\{\partial_\xi + \sumposroots k_\alpha \frac{(\alpha,\xi)}{(\alpha,\,.\,)} (1-\refl)\right\}\phi_\infty(\lambda,\,.\,)=(\xi,\lambda)\phi_\infty(\lambda_,\,.\,).
\]
This is the defining equation for the Dunkl kernel, which has a local solution
space spanned by $\Exp(\lambda,k,\,.\,)$. Since obviously
$\phi_\infty(\lambda,0)=1$, we conclude that
$\phi_\infty(\lambda,\,z)=\Exp(\lambda,k,z)$ for all $(\lambda,z)\in V_1\times
V_2$. But this implies that the $\phi_n$ converge uniformly to
$\Exp(\,.\,,k,\,.\,)$ on the compact set $S\subset V_1\times V_2$ after all,
which contradicts \eqref{eq:nonconvergence}.
\end{proof}
In the end, we will only need the special case where $\lambda\in i\ga$ of the
following corollary. Note that the function involved is defined for all
$\epsilon$ under consideration.

\begin{corollary}\label{cor:limittransition}
Let $k\geq 0$. Then, as $\epsilon\rightarrow 0$ through the real numbers, we
have $\lim_{\epsilon\rightarrow 0} C_G(\epsilon^{-1}\lambda,k,\epsilon x)=\Exp
(\lambda,k,x)$ for all $\lambda\in\gac$ and all $x\in\ga$.
\end{corollary}

We will now invoke Opdam's Paley--Wiener theorem. To this end, fix $x\in\ga$,
let $\textup{co}(G\cdot x)$ denote the convex hull of the orbit of $x$, and fix
$y\notin \textup{co}(G\cdot x)$. Suppose that $f\in {\mathcal
H}_{\textup{co}(G\cdot x)}$. Then, if $k_\alpha>0$ for all $\alpha$, the
Paley--Wiener theorem for the Cherednik transform \cite[Theorem
8.6.(2)]{gradhecke} implies that
\[
\int_\ga f(\lambda) C_G(i\lambda,k,y)\prod_{a\in\rplus} \left( 1-\frac{k_\alpha (\alpha,\alpha)}{2i(\lambda,\alpha)}\right)\prod_{\alpha\in\rplus}\left\vert\frac{\G\left(\frac{2i(\lambda,\alpha)}{(\alpha,\alpha)}+k_\alpha\right)}{\G\left(\frac{2i(\lambda,\alpha)}{(\alpha,\alpha)}\right)}\right\vert^2\,d\lambda=0,
\]
where the integrand should be read as $0$ for all $\lambda\in\ga^{\sing}$.
Now note, for real $\epsilon\neq 0$, that the map $\lambda\mapsto f(\epsilon\lambda)$ is an element of ${\mathcal H}_{\textup{co}(G\cdot \epsilon x)}$. Since $\epsilon y\notin \textup{co}(G\cdot \epsilon x)$, we have
\[
\int_\ga f(\epsilon\lambda) C_G(i\lambda,k,\epsilon y)\prod_{a\in\rplus} \left( 1-\frac{k_\alpha (\alpha,\alpha)}{2i(\lambda,\alpha)}\right)\prod_{\alpha\in\rplus}\left\vert\frac{\G\left(\frac{2i(\lambda,\alpha)}{(\alpha,\alpha)}+k_\alpha\right)}{\G\left(\frac{2i(\lambda,\alpha)}{(\alpha,\alpha)}\right)}\right\vert^2\,d\lambda=0.
\]
After a change of variables it follows that
\begin{equation}\label{eq:zerointegral}
\int_\ga f(\lambda) C_G(i\epsilon^{-1}\lambda,k,\epsilon y)\prod_{a\in\rplus} \left( 1-\frac{\epsilon k_\alpha (\alpha,\alpha)}{2i(\lambda,\alpha)}\right) \prod_{\alpha\in\rplus}\left\vert\epsilon^{k_\alpha}\frac{\G\left(\frac{2i\epsilon^{-1}(\lambda,\alpha)}{(\alpha,\alpha)}+k_\alpha\right)}{\G\left(\frac{2i\epsilon^{-1}(\lambda,\alpha)}{(\alpha,\alpha)}\right)}\right\vert^2\,d\lambda=0,
\end{equation}
for all $\epsilon>0$. Now it follows from a well-known result for the Gamma function \cite[p.~151]{titchmarsh} that, for each $a\in\R$, there exist constants $M_1(a),M_2(a),\widetilde M_1(a), \widetilde M_2(a)\geq 0$ such that, for all $s\in\R$,
\begin{align}\label{eq:dominatedestimates1}
\left\vert\frac{\Gamma(is+a)}{\Gamma(is)}\right\vert&\leq M_1(a)\vert s\vert^a + M_2(a),\\
\left\vert\left(1-\frac{a}{is}\right)\frac{\Gamma(is+a)}{\Gamma(is)}\right\vert&\leq \widetilde M_1(a)\vert s\vert^a + \widetilde M_2(a).\label{eq:dominatedestimates2}
\end{align}
Furthermore, for all $s,a\in\R,\,s\neq 0$, one has
\begin{equation}\label{eq:limitgamma}
\lim_{\epsilon\downarrow 0}\left\vert\epsilon^{a}\frac{\G(i\epsilon^{-1}s+a)}{\G(i\epsilon^{-1}s)}\right\vert=\vert s\vert^a.
\end{equation}
From \eqref{eq:dominatedestimates1} and \eqref{eq:dominatedestimates2} one sees
that the dominated convergence theorem applies in \eqref{eq:zerointegral} as
$\epsilon\downarrow 0$, since $f\in\schw(\ga)$ and $\vert
C_G(i\epsilon^{-1}\lambda,k,\epsilon y)\vert\leq\sqrt{\vert G\vert}$ by
\eqref{eq:cherednikestimate}. In the limit one thus obtains from
\eqref{eq:limitgamma} and Corollary~\ref{cor:limittransition} that
\begin{equation}\label{eq:expintegralzero}
\inta f(\lambda) \Exp( i\lambda,k,y)w_k(\lambda)\,d\lambda=0,
\end{equation}
for $k>0$, $x\in\ga$, $f\in{\mathcal H}_{\textup{co}(G\cdot x)}$, and $y\notin
\textup{co}(G\cdot x)$.

It is now an easy matter to prove the following.

\begin{theorem}[Paley--Wiener theorem, third version]\label{thm:pwhecke}
Let $G$ be a Weyl group and suppose that $\Re k\geq 0$. If $S\subset\ga$ is an
intersection of convex hulls of $G$-orbits, then the Dunkl transform is a
linear isomorphism between $\test(S)$ and ${\mathcal H}_{S}$.
\end{theorem}

\begin{proof}
As in the proof of Theorem~\ref{thm:convexpw}, in view of the inversion theorem
all that needs to be done is to prove that $E_k: {\mathcal
H}_{S}\mapsto\schw(\ga)$ has its image contained in $\test(S)$. To this end,
first assume that $S$ is the convex hull of one orbit and that $k>0$. In that
case, \eqref{eq:expintegralzero} is just the required result. But then
\eqref{eq:expintegralzero} actually holds for $\Re k>0$ by analytic
continuation in $k$, and the case $\Re k\geq 0$ follows from this again by
continuity. The result for an intersection of convex hulls of orbits follows
trivially from the result for the convex hull of one orbit.
\end{proof}

\begin{remark}
The type of limit transition above has also been introduced by Ben Sa\"id and \O rsted \cite{saidoersted}. The approach in [loc.cit.] starts from the invariant case,
using shift operators, and then proceeds to the general case. The use of shift
operators thus restricts the validity of the proofs to integral multiplicities.
The approach to this transition in the present paper appears to be somewhat
simpler, since Corollary~\ref{cor:limittransition} is seen to hold directly in
the general case and for more general multiplicities.
\end{remark}

\section{Extending the intertwining operator}\label{sec:intertwining}

In this section we return to the general situation of an arbitrary reflection
group and $\Re k\geq 0$. With the aid of the Paley--Wiener theorems, we will
define a linear isomorphism $V_k:C^\infty(\ga)\mapsto C^\infty(\ga)$, such that
$T_\xi(k)V_k=V_k\partial_\xi$ $(\xi\in\ga)$, and which extends the intertwiner
operator as it has originally been constructed for polynomials by Dunkl for $k\geq 0$
\cite{operatorscommuting}. The polynomial case for general regular $k$ was considered in \cite{singpol} and an extension of the intertwiner operator to an algebra of real
analytic functions was constructed in \cite{integralkernels} for $k\geq 0$.
The extension to $C^\infty(\ga)$ for $\Re k\geq 0$ as in this section was first established in \cite{thesis}. The case where $k\geq 0$ and $\sumposroots k_\alpha>0$ was later on also considered in \cite{trimeche1}, using R\"osler's representing measures which are presently only known to exist when $k\geq 0$.

Recall that the symmetry $S$, defined by $Sf(x)=f(-x)$, relates $\ek$ and $\dk$
by $S\dk=\dk S=\ek$ \cite[Lemma 4.3]{dunkltransform}. The Paley--Wiener
theorem~\ref{thm:pw} therefore implies that $\dk E_0=S\ek E_0$ maps
$\test(\ga)$ into itself. This map is easily seen to be continuous in the standard inductive limit topology (e.g., as a
consequence of the closed graph theorem) and it is actually a homeomorphism
since the inverse is $D_0 \ek$. Thus the transpose $(D_k E_0)^t:
\test^\prime(\ga)\mapsto\test^\prime(\ga)$ is a homeomorphism in the $\tu{weak}^*$-topology, leaving $\schw^\prime(\ga)$ invariant as a consequence of part (5) in
Theorem~\ref{thm:inventionesresults}.

For the formulation of the following theorem, we let $\openball$ denote the open ball of radius $R$, centered at the
origin.

\begin{theorem}\label{thm:intertwiner} Let $G$ be a finite reflection group and suppose that $\Re k\geq 0$. The map $W_k:C^\infty(\ga)
\mapsto\test^\prime(\ga)$, defined by
\[W_kf=\frac{c_0}{c_k}(\dk E_0)^t (fw_k),\]
is actually a
linear automorphism of $C^\infty(\ga)$. Let $V_k:C^\infty(\ga)\mapsto
C^\infty(\ga)$ be its inverse. Then both $W_k$ and $V_k$ commute with the
$G$-action and the following hold:
\begin{enumerate}
\item $W_k\tk=\partial_\x W_k\,\,(\x\in\ga)$. \item $W_k1=1$.
\item $W_k\poly_n=\poly_n$.
\item
\begin{enumerate}
\item If $f\in C^\infty(\ga)$ vanishes on $\openball$, then $W_k f$ vanishes on $\openball$.
\item If $U$ is an open $G$-invariant convex subset of $\ga$, and $f\in C^\infty(\ga)$ vanishes on $U$, then $V_k f$ vanishes on $U$.
\end{enumerate}
\item If $f\in C^\infty(\ga)\cap\lone$, $\dk
f\in\lone$, and $x\in\ga$, then
\[W_k f(x)=\frac{c_0}{c_k} E_0(w_k D_k f)(x)=\frac{1}{c_k}\inta \dk f(\lambda) \exp
(i(\lambda,x))w_k(\lambda)\,d\lambda.
\]
\item If $f\in
C^\infty(\ga)\cap L_1(\ga,dx)$, $D_0 f\in L_1(\ga,dx)$, and $x\in\ga$, then
\[
V_k f(x)=\frac{c_k}{c_0}E_k(w_k^{-1}D_0 f)=\frac{1}{c_0}\inta D_0f(\lambda)\Exp (ix,k,\lambda)\,d\lambda.
\]
\end{enumerate}
\end{theorem}

\begin{proof}
As a preparation, we recall from \cite[Lemma 4.13]{dunkltransform} that $\inta (\dk f) g w_k\,dx=
\inta f (\dk g)w_k\,dx$ for $f,g\in\lone$, and similarly for $\ek$.

Suppose $T\in\dist(\ga)$ and ${\supp}\,T\cap\openball=\emptyset$. Then
${\supp}\,(\dk E_0)^t T\cap\openball=\emptyset$. Indeed, let $\psi\in
\test(\openball)$. Then $\langle (\dk E_0)^t T,\psi\rangle=
\langle T,\dk E_0\psi\rangle=\langle T,\ek D_0\psi\rangle=0$ by the Paley--Wiener theorem~\ref{thm:pw}.

Now let $f\in
C^\infty(\ga)$. Fix $R>0$ and choose $\phi\in\test(\ga)$ such that $\phi=1$
on $\openball$. Then $\supp\, (\phi fw_k-fw_k)\cap\openball=\emptyset$,
hence by the above we have for arbitrary $\psi\in\test(\openball)$:
\begin{align}\label{eq:firstintegral}
\langle(\dk E_0)^t(fw_k),\psi\rangle&=\langle(\dk E_0)^t(\phi fw_k),\psi\rangle\\
&=\inta \psi(x) E_0\left(\dk(\phi f)w_k\right) (x)\,dx.\notag
\end{align}
Note that $D_k(\phi f)\in\schw(\ga)$, so $E_0\left(\dk(\phi f)w_k\right)$ is
smooth. Since the left hand side in \eqref{eq:firstintegral} does not depend on
$\phi$, the restriction of $E_0\left(\dk(\phi f)w_k\right)$ to $\openball$ is
apparently also independent of the choice of $\phi$. Hence we can unambiguously define a smooth function representing the distribution $(\dk E_0)^t (f w_k)$, by choosing for $x\in\ga$ any $\openball$ such that $x\in\openball$, any $\phi\in\test(\ga)$ such that $\phi=1$ on $\openball$, and putting $\left((\dk
E_0)^t (f w_k)\right)(x)= E_0\left(\dk(\phi f)w_k\right)(x)$. This
shows that $W_k$ maps $C^\infty(\ga)$ into itself. It follows from the remarks before the statement of the theorem that $W_k$ is injective.

In order to see that $W_k$ is also surjective
we consider $(D_0\ek)^t$. For $(D_0\ek)^t$ a geometrically more precise property with respect to supports can be proved than for $(D_kE_0)^t$. To be precise:  if $U$ is an open $G$-invariant convex subset of $\ga$, and $T\in\test^\prime(\ga)$ is such that $\supp\, T\cap U=\emptyset$, then $\supp\, (D_0\ek)^tT\cap U=\emptyset$. Indeed, if $\psi\in\test(U)$, then $\langle (D_0\ek)^t T,\psi\rangle=\langle T,D_0 \ek \psi\rangle=\langle T,E_0 \dk \psi\rangle=0$ by the geometrical form of the Paley--Wiener theorem for the ordinary Fourier transform---note that $\dk\psi$ is of Paley--Wiener type corresponding to the $G$-invariant compact convex set $\textup{co}(G\cdot\supp \,\psi)\subset U$. Let
$R>0$ and $f\in C^\infty(\ga)$ and choose $\phi$ as above. Taking the open ball $\openball$ as $U$, we then compute as follows for arbitrary $\psi\in\test(\openball)$:
\begin{align}\label{eq:secondintegral}
\langle(D_0 \ek)^t f,\psi\rangle&=\langle(D_0 \ek)^t(\phi f),\psi\rangle\\
&=\inta D_0(\phi f)(\lambda)(\ek\psi)(\lambda)\,d\lambda\notag\\
&=c_k^{-1} \inta D_0(\phi f)(\lambda)\left\{\inta \psi(x)\Exp (ix,k,\lambda) w_k(x)\,dx\right\}\,d\lambda\notag\\
&=\inta\left\{c_k^{-1}\inta D_0(\phi f)(\lambda)\Exp (ix,k,\lambda)\,d\lambda\right\}\psi(x)\,w_k(x)\,dx\notag
\end{align}

Now the function $g_\phi:x\mapsto c_k^{-1}\inta D_0(\phi f)(\lambda)\Exp
(ix,k,\lambda)\,d\lambda$ is smooth on $\ga$ since the derivatives of $\Exp
(i\,.\,,k,\lambda)$ are of polynomial growth in $\lambda$ \cite[Corollary
3.7]{dunkltransform} and $D_0(\phi f)\in\schw(\ga)$. The $g_\phi$ patch
together, as before, to a smooth function $g$ on $\ga$ which satisfies
$(D_0 E_k)^tf=gw_k$ by construction. Taking inverses we find $(\dk E_0)^t
(gw_k)=f$, showing that $W_k$ is surjective, as required.

It is easy to check that $W_k$ and hence its inverse $V_k$ commute with the $G$-action, using an invariant function $\phi$ in the construction, and noting that $D_k$ and $E_0$ have this property.

Continuing, we define the Dunkl operator $T_\x^\prime(k)$ on $\dist(\ga)$ by
$\langle T_\x^\prime(k) T,\psi\rangle=-\langle T,\tk \psi\rangle$. We have
the compatibility $T_\x^\prime(k)(fw_k)=(\tk f)w_k$ for $f\in C^\infty(\ga)$. Hence
\begin{align*}
\langle(\dk E_0)^t((\tk f)w_k),\psi\rangle
&=\langle (T_\x(k) f)w_k,\dk E_0\psi\rangle\\
&=\langle T_\x^\prime(k) (fw_k),\dk E_0\psi\rangle\\
&=-\langle fw_k,\dk E_0 \partial_\x\psi\rangle\\
&=\langle\partial_\x((\dk E_0)^t(f w_k)),\psi\rangle,\\
\end{align*}
proving part (1).

For part (2) and part (3) we consider $D_k^t:\tempdist(\ga)\mapsto\tempdist
(\ga)$. By the inversion theorem we have $D_k^t w_k=c_k\delta$, hence $(D_k E_0)^t
w_k=c_k E_0^t\delta=c_k c_0^{-1}$, proving part (2). Furthermore, one easily verifies the relation
$D_k^t (M_{\x^*} T)=iT_\x^\prime(k)D_k^t T$ for $T\in\tempdist(\ga)$ and $\x\in\ga$.
Thus for $p\in\poly_n$ we have $D_k^t (pw_k)=c_k i^n T_p^\prime(k)\delta$. By the
homogeneity of the $\tk$, there exists a differential operator $\widetilde{D}$
with constant coefficients and of homogeneous degree $n$ such that $c_k i^n
T_p^\prime(k)\delta=\widetilde{D}\delta$. Hence $(D_k E_0)^t (pw_k)=E_0^t\widetilde{D}\delta$ is in $\poly_n$, as stated in part (3).

Part (4) follows from the properties regarding the supports of distributions which we have used above.

The parts (5) and (6) follow by performing the computations in \eqref{eq:firstintegral} and \eqref{eq:secondintegral}
with $\phi=1$, which is validated by the integrability conditions.
\end{proof}

\begin{remark}\quad
\begin{enumerate}
\item
It is easily seen that the validity of Conjecture~\ref{conj:pw}
would imply that the open ball in part (4)(a) for $W_k$ could be replaced by an open $G$-invariant convex set, as for $V_k$ in part (4)(b). In fact, if all multiplicities are strictly
positive integers, then we already know that this stronger statement actually holds, in view of Theorem~\ref{thm:convexpw}. Likewise, if the groups is a Weyl group, then we known that this stronger statement holds for the interior of an intersection of convex hulls of orbits, in view of Theorem~\ref{thm:pwhecke}.
\item For $k\geq 0$, part (4)(b) also follows from the properties of R\"osler's representing measures, cf.~\cite{trimeche1}.
\end{enumerate}
\end{remark}

\section{Connection with the Cartan motion group}\label{sec:cmg}

If $G$ is a Weyl group $W$, then, for certain multiplicities, the $W$-invariant
part of the theory for Dunkl operators has an interpretation in terms of the
Cartan motion group, as we will explain in this section. The starting point is Heckman's result
\begin{equation}\label{eq:heckmanresult2}
T_p(k)=\frac{1}{n!}\left(\ad \,\frac{\D_k}{2}\right)^n M_p,
\end{equation}
for $p\in\poly_n$, cf.\ Section~\ref{sec:formularium}. A double application of this result enables us to
give an explicit description of the radial parts in terms of Dunkl operators.
From this description we can then conclude that the generalized Bessel functions coincide with the restriction of the spherical functions to $\ga$. This
identification, in turn, enables us to identify the restriction to the
invariants of the operators $W_k$ and $V_k$ from Theorem~\ref{thm:intertwiner} in terms of
(the flat analogue of) the Abel transform. Furthermore, we show that in
certain cases shift operators can be used to essentially invert the Abel
transform by an invariant differential operator.

Establishing terminology, let $\group$ be a connected non-compact semisimple
Lie group with finite center, maximal compact subgroup $\maxcompact$ and
corresponding Cartan decomposition ${\mathfrak g}={\mathfrak k}\oplus{\mathfrak
p}$ of the Lie algebra. The group $\group_0=\maxcompact\ltimes\gp$ acts on the
flat symmetric space $\group_0/\maxcompact\simeq\gp$ as isometries for the
Killing form. The group of isometries of $\gp$ thus obtained is known
as the Cartan motion group.

Choose a maximal abelian subspace $\ga\subset{\mathfrak p}$. Then there is a
restriction isomorphism $\Res_\ga^{\mathfrak p}: C^\infty({\mathfrak
p})^{\maxcompact}\mapsto C^\infty(\ga)^W$, with analogues for compactly supported smooth invariant functions and for invariant polynomials. Let $\Sigma$ be the restricted roots with multiplicities
$m_\alpha\,\,(\alpha\in\Sigma)$. We consider $\Sigma$ to be a subset of $\ga$
by means of the Killing form $(\,.\,,\,.\,)$ and we let $W$ be the Weyl group,
acting on $\ga$.

The spherical function $\psi_\lambda\,\,(\lambda\in\gac)$ on the symmetric
space $\gp$ satisfies the equations
$\partial(p)\psi_\lambda=p(\lambda)\psi_\lambda$ for all $p$ in the algebra
$S({\mathfrak p})^{\maxcompact}$ of $\maxcompact$-invariant polynomials on
$\mathfrak p$, where as usual $\partial(p)$ is the constant coefficient
differential operator corresponding to $p$. Furthermore, $\psi_\lambda(0)= 1$.
By \cite[Corollary 2.3, p.~402]{helref}, $\psi_\lambda$ is the \emph{unique}
$\maxcompact$-invariant function on ${\mathfrak p}$ with these two properties.

We let $R\subset\ga$ be the normalized root system of the finite reflection
group $W$ (so some roots in a component $BC_n$ of $\Sigma$ will coincide in
$R$). Choose the multiplicity function $k:R\mapsto\R$ as
$k_\alpha=1/4\sum_{\beta\in\R\alpha\cap\Sigma} m_\beta$. With these
multiplicities, we define the corresponding operators $T_p(k)$ (acting on
functions on $\ga$) for the finite reflection group $W$. The crucial
observation, with $\Rad$ denoting radial part, is the following:
\begin{equation}\label{eq:radialpartdescription}
\Rad (\partial(p))f=T_{\restoa p}(k)f\quad\left(p\in S({\mathfrak p})^{\maxcompact},\,\,
f\in C^\infty(\ga)^W\right).
\end{equation}

Before proving this relation, let us note that it shows how the algebra of
radial parts can be obtained from Dunkl operators, as follows. One computes the operator $T_q(k)$ for each polynomial $q$ in a set of fundamental invariants, and then restricts
it to the invariant functions on $\ga$. This restriction acts  as a
differential operator (a moment's thought will make this clear), and the differential operators which are thus obtained generate the algebra of radial
parts.

Turning to the proof of $\ref{eq:radialpartdescription}$, we first recall \cite{diffdiff} that
\[
\D_k=\D_0+2\sumposroots k_\alpha M_{(\alpha^*)^{-1}}\left\{\partial_\alpha-M_{(\alpha^*)^{-1}}
(1-\refl)\right\}.
\]
Let $\D_{\mathfrak p}$ denote the ordinary Laplacian on ${\mathfrak p}$. By the choice of
$k$, the expression for $\Rad\,\D_{\mathfrak p}$ in \cite[Proposition 3.13,
p.~270]{helref}, and the above expression for $\D_k$, we have $\Rad\,\D_{\mathfrak p}f=
\D_{k} f$ for $f\in C^\infty(\ga)^W$.

Now let $p\in S({\mathfrak p})^{\maxcompact}$ be homogeneous of
degree $n$ and suppose $f\in C^\infty(\ga)^W$. We apply \eqref{eq:heckmanresult2} first on $\mathfrak p$
(with $k=0$) and then on $\ga$ (for our particular choice of $k$) to prove \eqref{eq:radialpartdescription}:
\begin{align*}
\Rad (\partial(p))f&=\frac{1}{n!} \Rad \left(\left(\ad \,\frac{\D_{\mathfrak p}}{2}\right)^n M_p\right)f\\
&=\frac{1}{n!}\left\{\left(\ad \,\frac{\Rad \,\D_{\mathfrak p}}{2}\right)^n (\Rad\, M_p)\right\} f\\
&=\frac{1}{n!}\left\{\left(\ad \,\frac{\D_{k}}{2}\right)^n (M_{\restoa p})\right\} f\\
&=T_{\restoa p}(k)f.
\end{align*}

After having identified the radial parts in terms of Dunkl operators, it is now
easy to show how the Bessel kernel $J_W$, i.e., the $W$-invariant component
of ${\textup{Exp}}_W$ as in \eqref{eq:besselkerneldef}, is related to the restrictions of the spherical functions to
$\ga$. Indeed, if $\lambda\in\gac$, then we note that $J_W(\lambda,k,\,.\,)$ has
an extension $J_W^\ext(\lambda,k,\,.\,)$ to a smooth $\maxcompact$-invariant
function on ${\mathfrak p}$. This extension is evidently equal to $1$ in $0$,
and satisfies the equations $\partial(p)
J_W^\ext(\lambda,k,\,.\,)=p(\lambda)J_W^\ext(\lambda,k,\,.\,)$ for all $p\in
S({\mathfrak p})^{\maxcompact}$, as a consequence of
\eqref{eq:radialpartdescription} and the very definition of radial parts. By
the uniqueness mentioned above, we conclude that
$J_W^\ext(\lambda,k,x)=\psi_\lambda(x)$ for $\lambda\in\gac$ and $x\in\ga$.

In this context of the Cartan motion group, we will now proceed to identify the restrictions to the invariants $C^\infty(\ga)^W$ of the linear automorphisms $W_k,V_k:C^\infty(\ga)\mapsto C^\infty(\ga)$ from Theorem~\ref{thm:intertwiner}.
Recall \cite[p.~467]{helref} that
\begin{equation}\label{eq:intpasinta}
\intp f(x)\,dx=(2\p)^{\frac{\dim\,\gp}{2}}c_k^{-1}\inta \restoa f(x) w_{k}\,dx\quad(f\in C_c^\infty(\gp)^\maxcompact).
\end{equation}
The constant is determined by considering the Gaussian. Since we have identified the spherical functions as Bessel
functions, we see from this integral formula that the normalized spherical transform $\fourier$, defined by
\[
\fourier f(\lambda)=(2\p)^{-\frac{\dim\,\gp}{2}}\intp f(x)\psi_{-i\lambda}(x)\,dx\quad(\lambda\in\ga,\,f\in C_c^\infty(\gp)^\maxcompact),
\]
factors as $\fourier=D_{k}\restoa$. From this relation one sees that the
inversion, Plancherel and Paley--Wiener theorems for the spherical transform
follow from the corresponding theorems for $D_{k}$, when specialized to the
invariants (the first two admittedly being almost trivial for the Cartan motion
group).

Let $\gq$ be the orthoplement of $\ga$ in $\gp$ and define the (flat analogue of) the Abel transform:
\[
Af(x)=\intq f(x+q)\,dq\quad(f\in C_c^\infty(\gp)^\maxcompact,\,x\in\ga).
\]
Since the spherical transform can also be written as
\[
\fourier f(\lambda)=(2\p)^{-\frac{\dim\,\gp}{2}}\intp f(x)e^{-i(\lambda,x)}\,dx\quad(\lambda\in\ga,\, f\in C_c^\infty(\gp)^\maxcompact),
\]
Fubini's theorem yields an additional factorization of $\fourier$, namely $\fourier= (2\p)^{-\frac{\dim\,\gq}{2}}D_0A$. The
relation $(2\p)^{-\frac{\dim\,\gq}{2}}D_0 A f=D_k\restoa f$ for $f\in C_c^\infty(\gp)^\maxcompact$ and the
Paley--Wiener theorems for $D_0$ and $D_k$ together then show that the Abel transform
establishes a linear isomorphism between $C_c^\infty(\gp)^\maxcompact$ and
$C_c^\infty(\ga)^W$. It is easily verified that $A: C_c^\infty(\gp)^\maxcompact\mapsto C_c^\infty(\ga)^W$ is in fact a linear homeomorphism, so that the isomorphisms
$A^t:\test^\prime(\ga)^W\mapsto\test^\prime(\gp)^\maxcompact$ and
$(A^{-1})^t:\test^\prime(\gp)^\maxcompact\mapsto\test^\prime(\ga)^W$ are
defined.

Recall the symmetry $S$ from Section~\ref{sec:pw}, defined by $Sf(x)=Sf(-x)$,
and note that, as maps from $C_c^\infty(\ga)^W$ into itself, we have
\begin{align*}
D_{k}E_0&=D_{k}SD_0=(2\p)^{\frac{\dim\,\gq}{2}}D_{k}S\fourier A^{-1}
= (2\p)^{\frac{\dim\,\gq}{2}}D_{k}SD_{k}\restoa A^{-1}\\
&=(2\p)^{\frac{\dim\,\gq}{2}}D_{k}E_{k}\restoa A^{-1}
=(2\p)^{\frac{\dim\,\gq}{2}}\restoa A^{-1}.
\end{align*}
Let $\exttop:C^\infty(\ga)^W\mapsto C^\infty(\gp)^\maxcompact$ be the inverse of $\restoa$. Using \eqref{eq:intpasinta} we find, for $f\in C^\infty(\ga)^W$ and
$\psi\in C_c^\infty(\ga)^W$:
\begin{align*}
\langle W_{k} f,\psi\rangle_\ga &=(2\p)^{\frac{\dim\,\ga}{2}}c_k^{-1}\langle (D_{k}E_0)^t (fw_k),\psi\rangle_\ga\\
&=(2\p)^{\frac{\dim\,\ga}{2}}c_k^{-1}\langle fw_{k}, D_{k}E_0\psi\rangle_\ga\\
&=(2\p)^{\frac{\dim\,\gp}{2}}c_k^{-1}\langle f w_{k},\restoa A^{-1}\psi\rangle_\ga\\
&=\langle\exttop f, A^{-1}\psi\rangle_\gp\\
&=\langle (A^{-1})^t\exttop f,\psi\rangle_\ga.
\end{align*}
We conclude that $W_{k}=(A^{-1})^t\exttop$, as
maps with domain $C^\infty(\ga)^W$. For the right hand
side, we know a priori only that it maps $C^\infty(\ga)^W$ into $\test^\prime(\ga)^W$, but this equality and
Theorem~\ref{thm:intertwiner} show that it is actually a linear
automorphism of $C^\infty(\ga)^W$. It follows from this that
$(A^{-1})^t$, which a priori is only known to map
$C^\infty(\gp)^\maxcompact\subset\test^\prime(\gp)^\maxcompact$ into
$\test^\prime(\ga)^W$, establishes in fact a linear isomorphism between
$C^\infty(\gp)^\maxcompact$ and $C^\infty(\ga)^W$. This enables us to identify
the restriction of Dunkl's intertwiner to $C^\infty(\ga)^W$ as $V_{k}=
\restoa A^t$, where the fact that $A^t:C^\infty(\ga)^W\mapsto C^\infty(\gp)^\maxcompact$ is a linear isomorphism is then implicitly used in
the interpretation of the right hand side.

\begin{remark}For certain multiplicities of the restricted roots, the results on shift operators can be used to essentially invert the Abel transform, as follows. Assume that $k_\alpha=1/4\sum_{\beta\in\R\alpha\cap\Sigma} m_\beta$ is  an integer for all $\alpha\in R$.
Then \eqref{eq:complicatedformula} implies that there exists a $W$-invariant
differential operator $\widetilde{D}$ of order $1/4
\sum_{\alpha\in\Sigma}m_\alpha$, with poles along the singular hyperplanes only,
such that $E_{k}g=\widetilde{D}E_0g$ for all rapidly decreasing $g\in
C^\infty(\ga)^W$. Therefore, if $f\in C^\infty_c(\ga)^W$, we find that $\restoa
A^{-1}f=(2\p)^{-\frac{\dim\,\gq}{2}}D_{k} S D_0 f=(2\p)^{-\frac{\dim\,\gq}{2}}E_{k}D_0 f= (2\p)^{-\frac{\dim\,\gq}{2}}\widetilde{D} E_0 D_0 f= (2\p)^{-\frac{\dim\,\gq}{2}}\widetilde{D} f$.
This shows that the Abel transform is for such multiplicities essentially
inverted by the $W$-invariant differential operator $(2\p)^{-\frac{\dim\,\gq}{2}}\widetilde D$.
\end{remark}

\begin{remark}
In this section it has become clear that the theory of spherical functions for
the Cartan motion group is intimately connected with the invariant part of the
theory of Dunkl operators. It can, however, not be excluded that the general theory of these
operators is also relevant in this context, for conceptual reasons which have gone unnoticed
sofar. The author at least  finds it hard to believe that, e.g., the
description of the algebra of radial parts should not be the reflection of some more
direct connection. There are, in effect, two additional indications that such a
connection may exist.

The first of these is given by the estimates on the spherical functions
$\psi_\lambda$. One has the integral representation
\[
\psi_\lambda(x)=\int_{\maxcompact}e^{(\lambda, kx)}\,dk,
\]
for $\lambda\in\gac$and $x\in\ga$. The Kostant convexity theorem \cite[Theorem 10.2, p.~473]{helref} then
yields the estimate $\vert\psi_\lambda(x)\vert\leq\exp (\max_{w\in
W}\Re(w\lambda,x))$. But alternatively, this also follows from the estimate for $\textup{Exp}_W$ in Section~\ref{sec:notations_previous}
and the identification of the spherical functions as generalized Bessel functions, as explained in this section. This natural occurrence of $\textup{Exp}_W$ points to relevance
of the original Dunkl operators in the Cartan motion group case.

The second indication consists of Torossian's application of Dunkl operators in
the inversion of the Chevalley restriction isomorphism between the
$\maxcompact$-invariant polynomials on ${\mathfrak p}$ and the $W$-invariant
polynomials on $\ga$ \cite{torossian}. Here the original Dunkl operators make an appearance, not just the operators corresponding to invariant polynomials.
\end{remark}

\begin{remark}
When the results in this section are combined with similar results for the
Cherednik operators in Section~\ref{subsec:gradhecke} and with
Theorem~\ref{thm:limittransition}, one obtains a limit transition from the
spherical functions for a Riemannian symmetric space of the non-compact type to
the spherical functions for the corresponding Cartan motion group. This has
previously been proved in a direct fashion by Ben Sa\"id and \O rsted in
\cite{saidoersted}. This limit transition enables one to obtain explicit formulas for the
spherical functions for the Cartan motion group.
\end{remark}





\begin{thebibliography}{99}

\bibitem{abramowitz}
M.~Abramowitz and I.A.~Stegun, eds., \emph{Handbook of mathematical functions}.
Dover Publications, New York, 1970.


\bibitem{saidoersted}
S.~Ben Sa\"id and B.~\O rsted,
\emph{Bessel functions for root systems via the trigonometric setting}, Int.\ Math.\ Res.\ Not.\ {\bf 2005} (9), 551-585.


\bibitem{cherednik}
I.~Cherednik, \emph{A unification of Knizhnik--Zamolodchikov equations and
Dunkl operators via affine Hecke algebras}, Invent.\ Math.\ {\bf 106} (1991),
411-432.

\bibitem{diejenvinet}
J.F.~van Diejen and L.~Vinet, \emph{Calogero--Moser--Sutherland Models}. CRM
Series in Mathematical Physics, Springer, 2000.

\bibitem{reflectiongroups}
C.F.~Dunkl, \emph{Reflection groups and orthogonal polynomials on the sphere},
Math.\ Z.\ {\bf 197} (1988), 33-60.

\bibitem{diffdiff}
C.F. Dunkl, \emph{Differential-difference operators associated to reflection
groups}, Trans.\ Amer.\ Math.\ Soc.\ {\bf 311}, (1989), 167-183.

\bibitem{operatorscommuting}
C.F.~Dunkl, \emph{Operators commuting with Coxeter group actions on
polynomials}, in ``Invariant theory and tableaux," Minneapolis, 1988, \emph{IMA
Vol.\ Math.\ Appl.\ }{\bf 19}, 107-177.

\bibitem{integralkernels}
C.F.~Dunkl, \emph{Integral kernels with reflection group invariance}, Canad.\
J.\ Math.\ {\bf 43} (1991), 1213-1227.

\bibitem{hankeltransform}
C.F.~Dunkl, \emph{Hankel transforms associated to finite reflection groups}, in
``Proceedings of the special session on hypergeometric functions on domains of
positivity, Jack polynomials and applications," Tampa, 1991, \emph{Contemp.\
Math.\ }{\bf 138}, 123-138.


\bibitem{dunkla2}
C.F.~Dunkl, \emph{Intertwining operators associated to the group $S_3$}, Trans.\
Amer.\ Math.\ Soc.\ {\bf 347} (1995), 3347-3374.

\bibitem{dunklsymmetric1}
C.F.~Dunkl, \emph{Singular polynomials for the symmetric groups}, Int.\ Math.\ Res.\ Not.\ {\bf 2004} (67), 3607-3635.

\bibitem{dunklsymmetric2}
C.F.~Dunkl, \emph{Singular polynomials and modules for the symmetric groups}, Preprint (2005). ArXiv: math.RT/0501494.

\bibitem{singpol}
C.F.~Dunkl, M.F.E.~de Jeu~and E.M.~Opdam, \emph{Singular polynomials for finite
reflection groups},  Trans.\ Amer.\ Math.\ Soc.\ {\bf 346} (1994), 237-256.

\bibitem{dunklopdam}
C.F.~Dunkl and E.M.~Opdam, \emph{Dunkl operators for complex reflection
groups}, Proc.\ London Math.\ Soc.\ (3) {\bf 86} (2003), 70-108.


\bibitem{dunklxu}
C.F.~Dunkl and Y.~Xu, \emph{Orthogonal polynomials of several variables}.
Cambridge Univ.\ Press, 2001.

\bibitem{heckmanremark}
G.J.~Heckman, \emph{A remark on the Dunkl differential-difference operators},
in ``Harmonic analysis on reductive groups," Brunswick, 1989, \emph{Progr.\
Math.\ }{\bf 101}, 181-191.

\bibitem{heckmanbourbaki}
G.J.~Heckman, \emph{Dunkl operators}, S\'eminaire Bourbaki {\bf 828}, 1996-97,
Ast\'erisque {\bf 245} (1997), 223-246.

\bibitem{helref}
S.~Helgason, \emph{Groups and geometric analysis}. Academic Press, New York,
1984.

\bibitem{helref2}
S.~Helgason, \emph{Geometric analysis on symmetric spaces}. American
Mathematical Society, Providence, 1994.

\bibitem{hormander}
L.~H\"ormander, \emph{The analysis of linear partial differential operators I}.
Springer, 1983.

\bibitem{dunkltransform}
M.F.E.~de~Jeu, \emph{The Dunkl transform}, Invent.\ math.\ {\bf 113} (1993),
147-162.

\bibitem{thesis}
M.F.E.~de~Jeu, \emph{Dunkl operators}. Thesis, Leiden University, 1994.

\bibitem{densesubspaces}
M.F.E.~de~Jeu, \emph{Subspaces with equal closure}, Constr.\ Approx.\ {\bf 20}
(2004), 93-157.

\bibitem{lusztig}
G.~Lusztig, \emph{Affine Hecke algebras and their graded version}, J.\ Amer.\
Math.\ Soc.\ {\bf 2} (1989), 599-685.

\bibitem{someapplications}
E.M.~Opdam, \emph{Some applications of hypergeometric shift operators},
Invent.\ Math.\ {\bf 98} (1989), 1-18.

\bibitem{besseletc}
E.M.~Opdam, \emph{Dunkl operators, Bessel functions and the discriminant of a
finite Coxeter group}, Comp.\ Math.\ {\bf 85} (1993), 333-373.

\bibitem{gradhecke}
E.M.~Opdam, \emph{Harmonic analysis for certain representations of graded Hecke
algebras}, Acta Math.\ {\bf 175} (1995), 75-121.

\bibitem{opdamjapan}
E.M.~Opdam, \emph{Lecture notes on Dunkl operators for real and complex
reflection groups}, MSJ Memoirs {\bf 8}, Math.\ Soc.\ of Japan, Tokyo, 2000.

\bibitem{roeslerpositivity}
M.~R\"osler, \emph{Positivity of Dunkl's intertwining operator}, Duke Math.\
J.\ {\bf 98} (1999), 445-463.

\bibitem{roeslersummerschool}
M.~R\"osler, \emph{Dunkl operators: theory and applications}, in ``Orthogonal
polynomials and special functions," Leuven, 2002, \emph{Lecture Notes in Math.\
}{\bf 1817}, 93-135.

\bibitem{roeslerdejeu}
M.~R\"osler and M.F.E.~de~Jeu, \emph{Asymptotic analysis for the Dunkl kernel},
J.\ Approx.\ Theory {\bf 119} (2002), 110-126.

\bibitem{roeslervoitmarkov}
M.~R\"osler and M.~Voit, \emph{Markov processes related with Dunkl operators},
Adv.\ Appl.\ Math.\ {\bf 21} (1998), 575-643.

\bibitem{ruref}
W.~Rudin, \emph{Functional analysis}. Tata McGraw-Hill, New Delhi, 1973.

\bibitem{distributions}
L.~Schwartz, \emph{Th\'eorie des distributions}. Hermann, Paris, 1978.

\bibitem{thangaveluxu}
S.~Thangavelu and Y.~Xu, \emph{Convolution operator and maximal function for Dunkl transform}, Preprint (2004). ArXiv: math.CA/0403049 v3.

\bibitem{titchmarsh}
E.C.~Titchmarsh, \emph{The theory of functions}. Oxford Univer.\ Press, 1986.

\bibitem{torossian}
C.~Torossian, \emph{Une application des op\'erateurs de Dunkl au th\'eor\`eme
de restriction de Chevalley}, C.R.\ Acad.\ Sci.\ Paris {\bf 318} (1994),
895-898.

\bibitem{trimeche1}
K.~Trim\`eche, \emph{The Dunkl intertwining operator on spaces of functions and
distributions and integral representation of its dual}, Integral Transform.\
Spec.\ Funct.\ {\bf 12} (2001), 349-374.

\bibitem{trimeche2}
K.~Trim\`eche, \emph{Paley--Wiener theorems for the Dunkl transform and Dunkl
translation operators}, Integral Transform.\ Spec.\ Funct.\ {\bf 13} (2002),
17-38.

\bibitem{xu}
Y.~Xu, \emph{Integration of the intertwining operator for $h$-harmonic
polynomials associated to reflection groups}, Proc.\ Amer.\ Math.\ Soc.\ {\bf
125} (1997), 2963--2973.

\end{thebibliography}
\end{document}